\catcode`\^^Z=9
\catcode`\^^M=10
\output={\if N\header\headline={\hfill}\fi
\plainoutput\global\let\header=Y}
\magnification\magstep1
\tolerance = 500
\hsize=14.4true cm
\vsize=22.5true cm
\parindent=6true mm\overfullrule=2pt
\newcount\kapnum \kapnum=0
\newcount\parnum \parnum=0
\newcount\procnum \procnum=0
\newcount\nicknum \nicknum=1
\font\ninett=cmtt9

\font\ninebf=cmbx9

\font\sixbf=cmbx6
\font\ninesl=cmsl9

\font\nineit=cmti9

\font\ninerm=cmr9

\font\sixrm=cmr6
\font\ninei=cmmi9
\font\eighti=cmmi8
\font\sixi=cmmi6
\skewchar\ninei='177 \skewchar\eighti='177 \skewchar\sixi='177
\font\ninesy=cmsy9
\font\eightsy=cmsy8
\font\sixsy=cmsy6
\skewchar\ninesy='60 \skewchar\eightsy='60 \skewchar\sixsy='60
\font\titelfont=cmr10 scaled 1440
\font\paragratit=cmbx10 scaled 1200

\font\name=cmcsc10
\font\emph=cmbxti10

\font\tenmsbm=msbm10
\font\sevenmsbm=msbm7
%

%
\font\got=eufm10

\font\teneufm=eufm10
\font\seveneufm=eufm7
\font\fiveeufm=eufm5
\newfam\eufmfam
\textfont\eufmfam=\teneufm
\scriptfont\eufmfam=\seveneufm
\scriptscriptfont\eufmfam=\fiveeufm

\font\tenmsam=msam10
\font\sevenmsam=msam7
\font\fivemsam=msam5
\newfam\msamfam
\textfont\msamfam=\tenmsam
\scriptfont\msamfam=\sevenmsam
\scriptscriptfont\msamfam=\fivemsam
\font\tenmsbm=msbm10
\font\sevenmsbm=msbm7
\font\fivemsbm=msbm5
\newfam\msbmfam
\textfont\msbmfam=\tenmsbm
\scriptfont\msbmfam=\sevenmsbm
\scriptscriptfont\msbmfam=\fivemsbm
\def\Bbb#1{{\fam\msbmfam\relax#1}}
\def\cz{{\kern0.4pt\Bbb C\kern0.7pt}
}
\def\ez{{\kern0.4pt\Bbb E\kern0.7pt}
}
\def\fz{{\kern0.4pt\Bbb F\kern0.3pt}}
\def\gz{{\kern0.4pt\Bbb Z\kern0.7pt}}
\def\hz{{\kern0.4pt\Bbb H\kern0.7pt}
}
\def\kz{{\kern0.4pt\Bbb K\kern0.7pt}
}
\def\nz{{\kern0.4pt\Bbb N\kern0.7pt}
}
\def\oz{{\kern0.4pt\Bbb O\kern0.7pt}
}
\def\rz{{\kern0.4pt\Bbb R\kern0.7pt}
}
\def\sz{{\kern0.4pt\Bbb S\kern0.7pt}
}
\def\pz{{\kern0.4pt\Bbb P\kern0.7pt}
}
\def\qz{{\kern0.4pt\Bbb Q\kern0.7pt}
}
\newskip\ttglue
\def\ninepoint{\def\rm{\fam0\ninerm}%
  \textfont0=\ninerm \scriptfont0=\sixrm \scriptscriptfont0=\fiverm
  \textfont1=\ninei \scriptfont1=\sixi \scriptscriptfont1=\fivei
  \textfont2=\ninesy \scriptfont2=\sixsy \scriptscriptfont2=\fivesy
  \textfont3=\tenex \scriptfont3=\tenex \scriptscriptfont3=\tenex
  \def\it{\fam\itfam\nineit}%
  \textfont\itfam=\nineit
  \def\sl{\fam\slfam\ninesl}%
  \textfont\slfam=\ninesl
  \def\bf{\fam\bffam\ninebf}%
  \textfont\bffam=\ninebf \scriptfont\bffam=\sixbf
   \scriptscriptfont\bffam=\fivebf
  \def\tt{\fam\ttfam\ninett}%
  \textfont\ttfam=\ninett
  \tt \ttglue=.5em plus.25em minus.15em
  \normalbaselineskip=11pt
  \font\name=cmcsc9
  \let\sc=\sevenrm
  \let\big=\ninebig
  \setbox\strutbox=\hbox{\vrule height8pt depth3pt width0pt}%
  \normalbaselines\rm
  \def\sl{\it}}

\headline={\ifodd\pageno\rightheadline\else\leftheadline\fi}
\def\rightheadline{\ninepoint Paragraphen"uberschrift\hfill\folio}
\def\leftheadline{\ninepoint\folio\hfill Chapter"uberschrift}
\let\header=Y
\def\titel#1{\need 9cm \vskip 2truecm
\parnum=0\global\advance \kapnum by 1
{\baselineskip=16pt\lineskip=16pt\rightskip0pt
plus4em\spaceskip.3333em\xspaceskip.5em\pretolerance=10000\noindent
\titelfont Chapter \uppercase\expandafter{\romannumeral\kapnum}.
#1\vskip2true cm}\def\leftheadline{\ninepoint
\folio\hfill Chapter \uppercase\expandafter{\romannumeral\kapnum}.
#1}\let\header=N
}
\def\Titel#1{\need 9cm \vskip 2truecm
\global\advance \kapnum by 1
{\baselineskip=16pt\lineskip=16pt\rightskip0pt
plus4em\spaceskip.3333em\xspaceskip.5em\pretolerance=10000\noindent
\titelfont\uppercase\expandafter{\romannumeral\kapnum}.
#1\vskip2true cm}\def\leftheadline{\ninepoint
\folio\hfill\uppercase\expandafter{\romannumeral\kapnum}.
#1}\let\header=N
}
\def\need#1cm {\par\dimen0=\pagetotal\ifdim\dimen0<\vsize
\global\advance\dimen0by#1 true cm
\ifdim\dimen0>\vsize\vfil\eject\noindent\fi\fi}
\def\neupara#1{\par\penalty-2000
\procnum=0\global\advance\parnum by 1
\vskip1cm\noindent{\paragratit \the\parnum. #1}%
\def\rightheadline{\ninepoint\S\the\parnum.\ #1\hfill \folio}%
\vskip 8mm\noindent}
\def\Proclaim #1 #2\finishproclaim {\bigbreak\noindent
{\bf#1\unskip{}. }{\it#2}\medbreak\noindent}
%
\gdef\proclaim #1 #2 #3\finishproclaim {\bigbreak\noindent%
\global\advance\procnum by 1
{%
{\relax\ifodd \nicknum
\hbox to 0pt{\vrule depth 0pt height0pt width\hsize
   \quad \ninett#3\hss}\else {}\fi}%
\bf\the\parnum.\the\procnum\ #1\unskip{}. }
{\it#2}
\medbreak\noindent}
\newcount\stunde \newcount\minute \newcount\hilfsvar
\def\uhrzeit{
    \stunde=\the\time \divide \stunde by 60
    \minute=\the\time
    \hilfsvar=\stunde \multiply \hilfsvar by 60
    \advance \minute by -\hilfsvar
    \ifnum\the\stunde<10
    \ifnum\the\minute<10
    0\the\stunde:0\the\minute~Uhr
    \else
    0\the\stunde:\the\minute~Uhr
    \fi
    \else
    \ifnum\the\minute<10
    \the\stunde:0\the\minute~Uhr
    \else
    \the\stunde:\the\minute~Uhr
    \fi
    \fi
    }
 \def\calB{{\cal B}}

 \def\calD{{\cal D}}
 
 \def\calH{{\cal H}}
\def\calI{{\cal I}} \def\calJ{{\cal J}}

 \def\calR{{\cal R}}

 \def\calX{{\cal X}}

\def\gota{\hbox{\got a}} \def\gotA{\hbox{\got A}}

\def\dim{\mathop{\rm dim}\nolimits}

\def\GL{\mathop{\rm GL}\nolimits}

\def\id{\mathop{\rm id}\nolimits}
\def\im{\mathop{\rm Im}\nolimits} \def\Im{\im}

\def\kernel{\mathop{\rm kernel}\nolimits}

\def\mod{\mathop{\rm mod}\nolimits}
\def\O{{\rm O}}
\def\U{{\rm U}}

\def\proj{\mathop{\rm proj}\nolimits}

\def\SL{\mathop{\rm SL}\nolimits}
\def\SO{\mathop{\rm SO}\nolimits}

\def\Sp{\mathop{\rm Sp}\nolimits}

\def\boxit#1{
  \vbox{\hrule\hbox{\vrule\kern6pt
  \vbox{\kern8pt#1\kern8pt}\kern6pt\vrule}\hrule}}
\def\Boxit#1{
  \vbox{\hrule\hbox{\vrule\kern2pt
  \vbox{\kern2pt#1\kern2pt}\kern2pt\vrule}\hrule}}

\def\smallni{\smallskip\noindent }
\def\medni{\medskip\noindent }
\def\bigni{\bigskip\noindent }
\def\Isom{\mathop{\;{\buildrel \sim\over\longrightarrow }\;}}
\def\lo{\longrightarrow}

\def\loma{\longmapsto}

\def\spitz#1{\langle#1\rangle}
\def\imag{{\rm i}}
\def\pii{\pi {\rm i}}

\def\set#1{\bigl\{\,#1\,\bigr\}}

\def\square{\hbox{\hbox to 0pt{$\sqcup$\hss}\hbox{$\sqcap$}}}
\def\qed{\ifmmode\square\else{\unskip\nobreak\hfil
\penalty50\hskip3em\null\nobreak\hfil\square
\parfillskip=0pt\finalhyphendemerits=0\endgraf}\fi}
\def\pn{\the\parnum.\the\procnum}
\def\downmapsto{{\buildrel
        {\vbox{\hbox{\hskip.2pt$\scriptstyle-$}}}
        \over{\raise7pt\vbox{\vskip-4pt\hbox{$\textstyle\downarrow$}}}}}
\def\LK{1.1}

\def\Bho{1.4}
\def\HF{1.5}

\def\DStar{2.2}
\def\Bal{2.3}
\def\Rmf{2.4}

\def\ST{3.1}
\def\VExH{3.2}
\def\RSf{3.3}

\def\RA{4.2}

\def\Fwe{5.4}

\def\TRe{5.8}

\def\LFb{6.2}
\def\LiLm{6.3}
\def\PYns{6.4}
\def\nsY{6.5}

\def\HS{6.7}

\def\ZdcX{7.7}

\def\Rab{7.10}
\def\Sab{7.11}

\def\Uab{7.13}

\input xy \xyoption{all}
\def\today{\number\day.\ifcase\month\or Januar\or Februar\or M\"arz
\or April\or Mai\or Juni\or Juli\or August\or September
\or Oktober\or November\or Dezember\fi
\space \number\year}
\def\uhrzeit{
    \stunde=\the\time \divide \stunde by 60
    \minute=\the\time
    \hilfsvar=\stunde \multiply \hilfsvar by 60
    \advance \minute by -\hilfsvar
    \ifnum\the\stunde<10
    \ifnum\the\minute<10
    0\the\stunde:0\the\minute~Uhr
    \else
    0\the\stunde:\the\minute~Uhr
    \fi
    \else
    \ifnum\the\minute<10
    \the\stunde:0\the\minute~Uhr
    \else
    \the\stunde:\the\minute~Uhr
    \fi
    \fi
    }
\nopagenumbers
\nicknum=0      
\def\RAND#1{\vskip0pt\noindent\hbox to 0mm{\hss\vtop to 0pt{%
  \raggedright\ninepoint\parindent=0pt%
  \baselineskip=1pt\hsize=2cm #1\vss}}}
 \noindent\centerline{\titelfont The modular variety of hyperelliptic
curves of genus three}
\medni
 \def\rad{{\rm rad}}
\def\Bl{{\rm Bl}}
\def\Biproj{{\rm biproj}}

\font\dun=cmss10
\font\sevendun=cmss7
\vskip 1.5cm
\centerline{
\vbox{\openup-3pt\noindent\hsize=6cm{\it    Eberhard Freitag
\hfil\break Mathematisches Institut
\hfil\break Im Neuenheimer Feld 288
\hfil\break D69120 Heidelberg
\hfil\hfil\break}{\sevenrm
FREITAG@\vskip-1mm\noindent
\null\quad MATHI.UNI-HEIDELBERG.DE}}
\vbox{\openup-3pt\noindent\hsize=6cm{\it
       Riccardo Salvati Manni
\hfil\break Dipartimento di Matematica
\hfil\break Piazzale Aldo Moro, 2
\hfil\break I-00185 Roma, Italy
}\hfil\break{\sevenrm
SALVATI@\vskip-1mm\noindent
\null\quad MAT.UNIROMA1.IT}}\hfill
}
\def\text{\hbox}

\def\leftheadline{\ninepoint\folio\hfill
Modular forms}%
\def\rightheadline{\ninepoint Introduction\hfill \folio}%
\headline={\ifodd\pageno\rightheadline\else\leftheadline\fi}
\bigni
\centerline{2007}
\vskip8mm\noindent
\let\header=N%
\def\imag{{\rm i}}%
%
{\paragratit Introduction}%
\medni The modular variety  of non singular and complete hyperelliptic
curves  with level-two structure of  genus $3$
is a 5-dimensional quasi projective variety
which admits several standard compactifications.
The first one comes from the period map, which
realizes this variety as a sub-variety of the
Siegel modular variety of level two and genus three $\calH_3/\Gamma_3[2]$.
We denote the hyperelliptic locus by $\calI_3[2]$ and its closure in the
Satake compactification by
$$\overline{\calI_3[2]}\subset \overline{\calH_3/\Gamma_3[2]}.$$
The hyperelliptic locus has 36 irreducible (isomorphic)
components, which correspond to the
36 even theta characteristics. We denote the component of the characteristic
$m$ by
$$X=\overline{\calI_3[2]}^m.$$
One of the purposes of this paper will be to describe the
equations of $X$ in a suitable projective embedding and its Hilbert function
(\TRe).
It will turn out that $X$ is normal.
\smallskip
Two further models use the fact that hyperelliptic curves of genus three can
be obtained as coverings of a projective line with $8$ branch points.
The level two structure of the curve corresponds to an ordering of these 8 points.
This  leads
to the configuration space of 8 ordered points on a line.
There are two  other important compactifications of this configuration space.
The first one comes from geometric invariant theory using  so-called
semistable degenerated
point configurations in $(P^1)^8$. We denote this {\dun GIT}-compactification by
$Y$.
The equations of this variety in a suitable projective embedding
(by so-called Specht polynomials) are known.
For a general approach to the {\dun GIT}-model $(P_1)^n//\SL(2)$
we refer to the paper [HM1], [HM2].
In the special case $n=8$
this variety also can by
identified with a Baily-Borel compactified ball-quotient [Ko2].
$$Y=\overline{\calB/\Gamma[1-\imag]}.$$
We will describe these results in some detail and
obtain new proofs including some finer results for them.
We will see for example that the graded algebra describing
$(P_1)^8//\SL(2)$ is a Gorenstein ring (\Bho) and that it describes the full
ring of modular forms of $\Gamma[1-\imag]$ (\Rmf). This part is essentially due
to Kondo [Ko2], but we can  avoid the
theory of $K_3$-surfaces which is an essential part of Kondo's approach.
\smallskip
The other compactification uses the fact that families of marked projective lines
$(P^1,x_1,\dots,x_8)$ can degenerate to stable marked curves $(C,x_1,\dots,x_8)$
of genus 0.
\smallskip
We use the standard notation $\bar M_{0,8}$ for this compactification.
It is known that $\bar M_{0,n}$ for arbitrary $n$
is a smooth and projective fine moduli scheme.
\smallskip
In [Ka] (see also [AL])
it has been proved that  the birational map
$$\bar M_{0,n}\lo (P_1)^n//\SL(2)$$
is everywhere regular.
It is also known that there is natural  everywhere regular map
(Torelli map)
$$\bar M_{0,n}\lo\overline{\calH_g/\Gamma_g[2]}\quad
\hbox{(Satake compactification)}$$
with $n=2g+2$.
Hence in the case $n=8$ we have a diagram
$$\xymatrix{
  &\bar M_{0,8}\ar[dl]\ar[dr]&\\Y\ar@{-->}[rr]&
 &X\;.}$$
 The horizontal arrow  is only birational but not everywhere regular.
\smallskip
In this paper we find another realization of this triangle which uses the fact
that there are graded algebras (closely related to algebras of modular forms)
$A,B$ such that
$$X=\proj(A),\qquad Y=\proj(B).$$
The rational map $\xymatrix{Y\ar@{-->}[r]& X\\}$ is induced by a
homomorphism of graded
algebras
$$A\lo B.$$
This homomorphism rests on the theory of Thomae (19th century), in which the
thetanullwerte of hyperelliptic curves have been computed.
We use Mumfords approach [Mu] for the construction of this
homomorphism.
Using the explicit equations for $A,B$ we can compute the base locus of $Y\to X$.
It turns out to be the union of 56 $P^3$ (\LFb). Blowing up the base locus we get
a projective variety $Y^*$ and a diagram
$$\xymatrix{&
  Y^*\ar[dl]\ar[dr]&\\Y\ar@{-->}[rr]&
 &X\;,}$$
where the diagonal arrows are everywhere regular.
But the model $Y^*$ is singular.
The space $Y$ contains 35 special points
(corresponding to non-stable points in the {\dun GIT}-model or to the cusps in the
ball-quotient model). Their  inverse images in $Y^*$
are 35 disjoint rational surfaces.
Blowing up them we get a dominant {\it smooth\/} model $\tilde Y\to Y^*$ and
especially
a diagram
$$
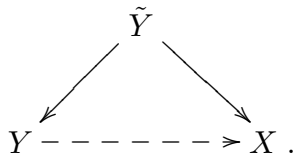
\xymatrix{&
  \tilde Y\ar[dl]\ar[dr]&\\Y\ar@{-->}[rr]&
 &X\;.}$$
We will see that $\tilde Y$ and $\bar M_{0,8}$ are isomorphic, \HS.
This can be considered as an explicit description of $\bar M_{0,8}$.
\smallskip
The variety $\tilde Y$ gives a correspondence between subvarieties of
  $Y=\bar X(8)^{\hbox{\sevendun GIT}}$ and of $Y=\overline{\calB/\Gamma[1-\imag]}$.
This correspondence explains several combinatorial similarities between the
two models.
These similarities can be described best, if one uses the ball-model to describe
$Y$. The reason is that in the description of the Siegel and the ball-model the
space $\fz_2^6$ occurs. In the Siegel case it occurs as a symplectic space,
its elements are so-called characteristics. In the ball case it occurs as quadratic
space. The symmetry of both models can be explained by the fact that in
characteristic 2 the symplectic and orthogonal world come together.
This will be described in the last section.
\smallskip
The second named author is grateful to Corrado De Concini and Shigeyuki Kondo for helpful  and enlightening discussions.
He would like also to thank Herbert Lange and Hanspeter Kraft for bringing  his attention to the papers
[AL] and [Ho].
\neupara{The configuration space for eight points
on a line}%
We recall some basic facts about the configuration space of eight
points on the line. With the help of computer computations we get very easy proofs
of them and also some new insight.
\smallskip
We consider the subset
$$\calX(8)\subset P^1(\cz)^8$$
of ordered eight tuples of elements of the projective
line, which are pairwise distinct.
The elements of $P^1(\cz)$ can be represented by columns
$a\choose b$, which are different from zero and
which are defined up to a constant factor.
Hence the elements of $\calX(8)$ can be represented
by matrices
$$M:=\pmatrix{a_1& a_2&\dots &a_8\cr b_1& b_2&\dots &b_8}.$$
Since $\SL(2,\cz)$ acts
on $P^1(\cz)$, we get an action of $\SL(2,\cz)$ on $\calX(8)$
by multiplication from the left.
The configuration space is defined as
$$X(8)=\SL(2,\cz)\backslash \calX(8).$$
We recall that a tableau is a matrix
$$\pmatrix{i_1\dots i_4\cr j_1\dots j_4}$$
which contains as entries all digits $i$, $1\le i\le 8$, and with the
property
$$i_1<i_2<i_3<i_4,\qquad i_1<j_1,\quad i_2<j_2,\quad i_3<j_3,\quad i_4<j_4.$$
There are 105 tableaux.
The tableau is called standard if in addition
$j_1<j_2<j_3<j_4$.
There are 14 standard tableaux:
$$\matrix{
\pmatrix{1&2&3&4\cr5&6&7&8}&
\pmatrix{1&2&3&5\cr4&6&7&8}&
\pmatrix{1&2&3&6\cr4&5&7&8}&
\pmatrix{1&2&3&7\cr4&5&6&8}\cr
\pmatrix{1&2&4&5\cr3&6&7&8}&
\pmatrix{1&2&4&6\cr3&5&7&8}&
\pmatrix{1&2&4&7\cr3&5&6&8}&
\pmatrix{1&2&5&6\cr3&4&7&8}\cr
\pmatrix{1&2&5&7\cr3&4&6&8}&
\pmatrix{1&3&4&5\cr2&6&7&8}&
\pmatrix{1&3&4&6\cr2&5&7&8}&
\pmatrix{1&3&4&7\cr2&5&6&8}\cr
\pmatrix{1&3&5&6\cr2&4&7&8}&
\pmatrix{1&3&5&7\cr2&4&6&8}\cr}$$
To each tableau one associates the expression
$$D(M)=\prod_{\nu=1}^4(a_{i_\nu}b_{j_\nu}-a_{j_\nu}b_{i_\nu}).$$
This leads to
the so-called Specht polynomial
$$D\pmatrix{1&\dots&1\cr X_1&\dots &X_8}=
(X_{i_1}-X_{j_1})(X_{i_2}-X_{j_2})(X_{i_3}-X_{j_3})(X_{i_4}-X_{j_4}).$$
The space generated by all the 105 Specht polynomials has dimension
14. A basis is given by the Specht polynomials of the standard
tableaux. We denote them by $Y_1,\dots,Y_{14}$
in the above ordering.

It is clear that these polynomials define a holomorphic map
$$X(8)\lo P^{13}(\cz).$$
It is easy to check and well-known that $X(8)\hookrightarrow P^{13}$ is
a smooth embedding.
We denote by $\bar X(8)$ the closure of the image of $X(8)$
in $P^{13}(\cz)$.
Hence
$$\bar X(8)=\proj \cz[Y_1,\dots,Y_{14}]/\cal J,$$
where $\calJ$ is the ideal of relations between the Specht polynomials.
The ideal $\calJ$ can be computed by elimination:
Consider the polynomial ring in 28 variables
$\cz[X_1,\dots,X_{14},T_1,\dots,T_{14}]$ and the ideal generated
by the relations $T_i-Y_i$  ($1\le i\le 14$). Then $\calJ$ is the
intersection of this ideal with $\cz[Y_1,\dots,Y_{14}]$. With computer
algebra one can verify now\footnote{*)}{The computer algebra calculation
has to done over the field $\qz$ instead $\cz$, which
is sufficient.} a slightly modified result of Koike [Koi]:
\proclaim
{Lemma (Koike)}
{The ideal $\calJ$ is generated by the polynomials
$$\eqalign{
&J_1=Y_{9}Y_{13}-Y_{8}Y_{14}\cr
&J_2=Y_{2}Y_{12}-Y_{3}Y_{12}-Y_{5}Y_{12}+Y_{6}Y_{12}-Y_{1}Y_{14}+Y_{3}Y_{14}+
Y_{5}Y_{14}-Y_{7}Y_{14}-\cr&\quad Y_{8}Y_{14}+Y_{9}Y_{14}-Y_{11}Y_{14}+Y_{13}Y_{14}\cr
&J_3=Y_{7}Y_{11}-Y_{6}Y_{12}\cr
&J_4=Y_{4}Y_{11}-Y_{3}Y_{12}-Y_{4}Y_{13}+Y_{3}Y_{14}\cr
&J_5=Y_{2}Y_{11}-Y_{5}Y_{11}-Y_{3}Y_{12}+Y_{6}Y_{12}-Y_{1}Y_{13}+
Y_{5}Y_{13}-Y_{7}Y_{13}+Y_{3}Y_{14}-\cr&\quad Y_{11}Y_{14}+Y_{13}Y_{14}\cr
&J_6=Y_{7}Y_{10}-Y_{9}Y_{10}-Y_{5}Y_{12}+Y_{5}Y_{14}\cr
&J_7=Y_{6}Y_{10}-Y_{8}Y_{10}-Y_{5}Y_{11}+Y_{8}Y_{11}-Y_{9}Y_{11}+Y_{5}Y_{13}-Y_{6}Y_{13}+Y_{6}Y_{14}\cr
&J_8=Y_{4}Y_{10}-Y_{9}Y_{10}-Y_{3}Y_{12}-Y_{5}Y_{12}+Y_{6}Y_{12}-Y_{4}Y_{13}-
Y_{1}Y_{14}+\cr&\quad Y_{3}Y_{14}+Y_{4}Y_{14}+Y_{5}Y_{14}-
Y_{7}Y_{14}+Y_{10}Y_{14}-Y_{11}Y_{14}+Y_{14^2}\cr
&J_9=Y_{3}Y_{10}-Y_{8}Y_{10}-Y_{5}Y_{11}+Y_{8}Y_{11}-Y_{9}Y_{11}-Y_{3}Y_{12}+
Y_{6}Y_{12}-Y_{1}Y_{13}+\cr&\quad Y_{5}Y_{13}- Y_{7}Y_{13}+Y_{10}Y_{13}-
Y_{11}Y_{13}+Y_{3}Y_{14}+Y_{13}Y_{14}\cr
&J_{10}=Y_{2}Y_{7}-Y_{3}Y_{7}-Y_{1}Y_{9}+Y_{3}Y_{9}-Y_{9}Y_{11}-
Y_{5}Y_{12}+Y_{6}Y_{12}+Y_{5}Y_{14}-\cr&\quad Y_{7}Y_{14}+Y_{9}Y_{14}\cr
&J_{11}=Y_{4}Y_{6}-Y_{3}Y_{7}-Y_{4}Y_{8}+Y_{7}Y_{8}+Y_{3}Y_{9}-Y_{6}Y_{9}-Y_{7}Y_{13}+Y_{6}Y_{14}\cr
&J_{12}=Y_{2}Y_{6}-Y_{3}Y_{7}-Y_{1}Y_{8}+Y_{3}Y_{9}-Y_{5}Y_{11}-
Y_{9}Y_{11}+Y_{6}Y_{12}+Y_{5}Y_{13}-\cr&\quad Y_{7}Y_{13}+Y_{8}Y_{14}\cr
&J_{13}=Y_{4}Y_{5}-Y_{3}Y_{7}-Y_{4}Y_{8}+Y_{7}Y_{8}-Y_{1}Y_{9}+Y_{3}Y_{9}+
Y_{4}Y_{9}-Y_{7}Y_{9}-Y_{9}Y_{11}-\cr&\quad Y_{5}Y_{12}+Y_{6}Y_{12}-
Y_{7}Y_{13}+Y_{5}Y_{14}+Y_{9}Y_{14}\cr
&J_{14}=Y_{3}Y_{5}-Y_{3}Y_{7}-Y_{1}Y_{8}+Y_{3}Y_{9}-Y_{5}Y_{11}-
Y_{9}Y_{11}+Y_{6}Y_{12}+Y_{5}Y_{13}-\cr&\quad Y_{7}Y_{13}+Y_{8}Y_{14}\cr
}$$
}
LK%
\finishproclaim
We want to describe the boundary $\bar X(8)-X(8)$.
Koike
 compares $\bar X(8)$ with a certain ball-quotient.
This relies to a paper of Kondo who identifies the Satake compactified
ball-quotient with the {\dun GIT}-compactification of $X(8)$.
Since Kondo's argument is not completely correct (s.~section 2),
we give here in some detail
a self contained complete description of $\bar X(8)$.
We need some results from geometric invariant theory which
gives the so-called {\dun GIT}-compactification of $X(8)$.
\smallskip
Recall that this is a certain projective variety $\proj R^G$. Here $R$
is the graded algebra $R=\sum_{n=0}^\infty R_n$, where
$R_n$ denotes the space of all polynomials in the entries of a
$2\times 8$-matrix $M$ which are homogenous of the same degree $n$
in each column of $M$. Hence $\proj(R)=P^1(\cz)^8$.
There is a natural action of the group $G=\SL(2,\cz)$ on $R$.
The ring of invariants $R^G$ is finitely generated. One defines
$$\bar X(8)^{\hbox{\sevendun GIT}}:=\proj(R^G).$$
Recall that a point from $P^1(\cz)^8$ is called semistable,
if there is some element from $R^G_n$, which doesn't vanish on it.
The set of all semistable points $P^1(\cz)^8_{\sevenrm ss}\subset
P^1(\cz)^8$ is a $G$-invariant Zariski open subset.
There is an obvious map $P^1(\cz)^8_{\sevenrm ss}\to \bar
X(8)^{\hbox{\sevendun GIT}}$
and this is a categorical quotient. Since we are in a geometric situation, this
map is surjective.
\smallskip
It is known that a point in $P^1(\cz)^8$ is semistable if and only
if not more then 4 of its components can agree.
From this description it follows that the Specht polynomials
don't have a joint zero. Hence we have a map
$P^1(\cz)^8_{\sevenrm ss}\to \bar X(8)$ and
from the universal property we get a map
$$\bar X(8)^{\hbox{\sevendun GIT}}\lo \bar X(8).$$
We will see that this is a biholomorphic map.
This actually follows from a known theorem from Kempe, s.~[HM1], [HM2]:
\smallni
{\it The ring $R^G$ is generated by its elements of degree one.
(These are the Specht polynomials).}
\smallni
We don't need to use this theorem. We will obtain is as consequence
of a sharper result which seems to be new (\Bho).
\smallskip
First of all one sees from the description of the semistable points
that the irreducible subvarieties of dimension 4
of the boundary $\bar X(8)^{\hbox{\sevendun GIT}}-X(8)$ are
defined by the subsets of $P^1(\cz)^8$ where two points agree. Hence these
irreducible subsets perform one orbit under the group $S_8$. It follows that
the boundary $\bar X(8)-X(8)$ also contains at most one orbit
of irreducible subvarieties of dimension 4. It is really one orbit, since
one can check easily by means of computer algebra that the subvariety
defined by $J_9=J_{10}=J_{11}=J_{12}=J_{13}=J_{14}=0$ is of dimension $4$.
It is easy to exhibit a special point in this subvariety which is smooth
in the whole $\bar X(8)$. Hence we see:
\proclaim
{Lemma}
{The ring $\cz[Y_1,\dots,Y_{14}]/\cal J$ is regular in codimension one.}
Rco%
\finishproclaim
Another calculation with computer algebra gives an explicit minimal resolution of
this ring:
\proclaim
{Lemma}
{Let $A$ be the polynomial ring $\cz[Y_1,\dots,Y_{14}]$
The minimal resolution of the ring $B:=\cz[Y_1,\dots,Y_{14}]/\cal J$
is of the form
$$\eqalign{
0\to A^{1}\to A^{14}&\to A^{175}\to A^{512}\to A^{700}\cr
\to &A^{512}\to A^{175}\to A^{14}\to A^1\to B\to 0.\cr}$$
{\bf Corollary.} The ring $B$ is a Cohen Macauley ring, even more, it
is a Gorenstein ring.}
CM%
\finishproclaim
Using Serre's criterion for normality we get now:
\proclaim
{Theorem}
{The ring $\cz[Y_1,\dots,Y_{14}]/\cal J$ is a normal Gorenstein ring.
\smallni
{\bf Corollary 1 \it(Kempe)\bf.} This ring agrees with the full ring of invariants
$R^G$ which hence is generated by the elements of lowest degree
(hence by the Specht polynomials).
\smallni
{\bf Corollary 2 \it(Kempe)\bf.} The map
$$\bar X(8)^{\hbox{\sevendun GIT}}\Isom \bar X(8)$$
is biholomorphic.}
Bho%
\finishproclaim
We want to point out that Kempe's result has been generalized in [HM1]
and [HM2] to point configurations of an arbitrary number of points
in the projective line.
\smallskip
From the computation of the resolution one can also read off:
\proclaim
{Theorem}
{The Hilbert function of $R^G\cong
\cz[Y_1,\dots,Y_{14}]/\cal J$ is
$$\eqalign{
\sum_{n=0}^\infty \dim R_n^G\; t^n&={1+8t+22t^2+8t^3+t^4\over
(1-t)^6}\cr
&=1+14t+91t^2+364t^3+1085t^4+5719t^6+\cdots\cr}$$}
HF%
\finishproclaim
This formula is in concordance with other formulae
in literature. Howe [Ho] for example gives
$$\dim R_n^G={n^5+5n^4+11n^3+13n^2+9n+3\over 3}\;,$$
which is obviously the same. The fact that Howe's formula is true for
all $n$, reflects that $R^G$ is Gorenstein.
\smallskip
Another remarkable formula has been communicated to us
by De Concini
$$\eqalign{
\dim R_n^G=
&\sum_{i_1+2i_2+\dots+ni_n=4n}{8 \choose i_1 }{8-i_1 \choose i_2
}{8-i_1-\dots-i_n \choose i_n }\cr&-
\sum_{i_1+2i_2+\dots+ni_n=4n-1}{8 \choose i_1 }{8-i_1 \choose i_2
}{8-i_1-\dots-i_n \choose i_n }\cr}$$

\neupara{A ball quotient}%
The ring $\cz[Y_1,\dots,Y_{14}]/\cal J$ of the previous section turns out to be
the ring of modular forms on a certain ball quotient. This is related to the fact
that $\bar X(8)$ is isomorphic to a Baily-Borel compactified ball-quotient.
These results
are essentially due to Kondo [Ko2]. Here we amend and extend some of his arguments.
\smallskip
The configuration space of 8 points on a line is related to the
even lattice
$$L=U\oplus U(2)\oplus D_4(-1)\oplus D_4(-1).$$
Here as usual $U=\gz^2$ with the quadratic form $x_1x_2$ and
$$D_4=\set{(x_1,x_2,x_3,x_4)\in\gz^4;\quad x_1+\cdots+x_4\equiv0\mod 2},\quad
q(x)=x_1^2+\cdots+x_4^2.$$
The notation $M(m)$ means  that one takes the same
abelian group $M$, but the quadratic form is multiplied by $m$.
Hence $\U\oplus \U(2)=\gz^4$ with the quadratic form $x_1x_2+2x_3x_4$.
The discriminant group $L'/L$ is isomorphic
to $\fz_2^6$. The isomorphism can be chosen such that the diagram
$$\matrix{ L'/L&\buildrel q\over\lo &\qz/\gz\cr
\uparrow&&\uparrow\cr
\fz_2^6&\lo&\fz_2\cr
(x_1,\dots,x_6)&\loma& x_1x_4+x_2x_5+x_3x_6\cr}$$
commutes. Here $\fz_2$ is embedded into $\qz/\gz$ by
sending $0\mapsto 0$ and $1\mapsto 1/2$.
It will be basic for us that the orthogonal
group (so-called even type) $\O(\fz_2^6)$ is isomorphic to the symmetric
group $S_8$.
We choose an isomorphism
$$S_8\Isom\O(\fz_2^6).$$
A natural concrete construction of such an isomorphism will be
given in section 4. At the moment the choice is not important.
\smallskip
We now describe a certain graded algebra, which is related to
$\O(\fz_2^6)$. This is the third member of a a sequence of algebras
related to $\O(\fz_2^{2m})$ [FS]. We recall very briefly its definition.
We attach to each maximal totally isotropic subspace $I\subset \fz_2^{2m}$
(of dimension $m$) a variable $X_I$ and we consider the subring
$$\calR_m:=\cz[\dots X_I-X_J\dots]$$
of the polynomial ring $\cz[\dots X_I\dots]$ in all
these variables. We define an ideal $\calI_m\subset \calR_m$.
It is the sum of a linear part and a quadratic part,
$$\calI_m=\calI_m^{\hbox{\sevenrm lin}}+\calI_m^{\hbox{\sevenrm qua}}.$$
We first define the linear part:
For this we need the characteristic function
$$\chi_I:\fz_2^{2m}\lo\cz$$
of a maximal totally isotropic subspace.
Let $A\subset\fz_2^{2m}$ be a totally
isotropic subspace of dimension $m-2$. There are 6 maximal totally isotropic
subspaces $I_1,\dots, I_6$ containing $A$. Their ordering can be chosen such that
$$(\chi_{I_1}-\chi_{I_2})+(\chi_{I_3}-\chi_{I_4})+(\chi_{I_5}-\chi_{I_6})=0.$$
The linear part $\calI_m^{\hbox{\sevenrm lin}}$ then is generated by all
$$(X_{I_1}-X_{I_2})+(X_{I_3}-X_{I_4})+(X_{I_5}-X_{I_6}).$$
Next we define the quadratic part: Again we consider
a totally
isotropic subspace $A\subset\fz_2^{2m}$ of dimension $m-2$ and the
6 maximal totally isotropic
subspaces $I_1,\dots, I_6$ containing it. After a suitable choice of their ordering
the relation
$$(\chi_{I_1}-\chi_{I_2})(\chi_{I_1}-\chi_{I_4})=
(\chi_{I_3}-\chi_{I_6})(\chi_{I_5}-\chi_{I_6})$$
holds. The quadratic part $\calI_m^{\hbox{\sevenrm qua}}$ of the ideal
is generated by all
$$(X_{I_1}-X_{I_2})(X_{I_1}-X_{I_4})-
(X_{I_3}-X_{I_6})(X_{I_5}-X_{I_6}).$$
The ring $\calR_m/\calI_m$ seems to be very interesting.
For example the cases $m=5$ and $m=6$ are related to Enriques surfaces.
Here we are interested in the case $m=3$.
\smallskip
First of all one computes that a linear combination $\sum_VC_VX_V$
is in the ideal $\calI_3$ if an only if  $\sum_VC_V\chi_V=0$.
The dimension of the
space generated by all $\chi_V-\chi_W$ can be computed and it is 14.
Hence $\calR_3/\calI_3^{\hbox{\sevenrm lin}}$ is a polynomial ring
in 14 variables.
THe group $\O(\fz_2^6)$ acts on this ring. The action comes from an irreducible
linear representation on the lowest degree part.
Since the subgroup of index two  $A_8\subset S_8$ admits
only one irreducible
representation of dimension 14,
we obtain an essentially unique $A_8$-equivariant
isomorphism
$$\calR_3/\calI_3^{\hbox{\sevenrm lin}}\lo\cz[Y_1,\dots,Y_{14}],$$
where on the right hand
side we have the ring of section one.
It turns out to be $\O(\fz_2^6)=S_8$-equivariant.
Now a concrete check shows that the quadratic relations go to 0.
Hence we get a homomorphism
$$\calR_3/\calI_3 \lo\cz[Y_1,\dots,Y_{14}]/\calJ.$$
It is no problem  to check by a concrete calculation:
\proclaim
{Theorem}
{The homomorphism
$$\calR_3/\calI_3 \lo\cz[Y_1,\dots,Y_{14}]/\calJ$$
is an $\O(\fz_2^6)\cong S_8$-equivariant isomorphism.
\smallni
{\bf Corollary.} The ring $\calR_3/\calI_3$ is a
normal ring of Krull dimension $5$.}
hiI%
\finishproclaim
We want to mention that the dimension of $\calR_m/\calI_m$
is unknown for $m>3$.
\smallskip
The ring $\calR_3/\calI_3$ is related to modular forms [FS]:
There is a subgroup of index two of $\O(L\otimes_\gz\rz)$, which doesn't contain
the reflections along vectors of po\-si\-tive norm. Its intersection with $\O(L)$
is denoted by $\O'(L)$. The dis\-cri\-mi\-nant kernel is
$$\Gamma_L=\kernel(\O'(L)\lo \O(L'/L)).$$
The Borcherds additive lifting attaches to each maximal totally isotropic subspace
a modular form [FS].
This gives a homomorphism
$$\calR_3/\calI_3^{\hbox{\sevenrm lin}}\lo A(\Gamma_L),$$
where $A(\Gamma_L)$ denotes the ring of modular forms with respect
to $\Gamma_L$. It is
possible to derive certain quadratic relations [FS]. Kondo derived
quartic relations but Koike [Koi] pointed out that they are
consequences of quadratic relations. These relations can be explained
naturally ba means of the
theory of Borcherds products: To explain this, we have to recall the
notion of a star [FS].
\proclaim
{Definition}
{A  star $S$ is a set of $4$ anistropic
vectors in $\fz^6$ which is a coset of a 2-dimensional totally
isotropic space.}
DStar%
\finishproclaim
(Kondo considers instead of stars ``maximal
totally singular subspaces''.  They are in one-to one- correspondence
with the stars. The set of anisotropic elements of such a space is a
star and this gives a one-to-one correspondence between stars and
maximal totally singular subspaces.) To every star a certain
Borcherds product can be associated. To explain this we have to
recall how to associate to an element $\alpha\in \fz_2^6$ of
non-zero norm a certain Heegner divisor $\calH_\alpha$ on the
symmetric domain which is associated to $\O(L)$.
Recall that the
elements of this symmetric domain can be considered as two
dimensional positive definite subspaces of $L\otimes_\gz\rz$. Then
 $\calH_\alpha$
consists of all such subspaces, which are orthogonal
to some
$$\delta\in L^*,\quad (\delta,\delta)=-1,\quad
\delta\ \hbox{over}\ \alpha.$$
If $S$ is  star, we denote by $\calH_S$ the
union of the four $\calH_\alpha$, $\alpha\in S$.
The basis fact is that the additive lift space contains for every star $S$ a form
$f_S$,
whose zero divisor contains $\calH_S$.
Since there is also a Borcherds product of the same weight with this property, the zero divisor equals $\calH_S$.
From this it is possible to derive directly
quadratic relations [FS] and to prove:
\proclaim
{Lemma}
{The Borcherds additive lifting gives a homomorphism
$$\calR_3/\calI_3\lo A(\Gamma_L).$$}
Bal%
\finishproclaim
Of course this homomorphism is not an isomorphism, since the left hand side
has Krull dimension 6 and the right hand side has Krull dimension 11.
The picture remedies if one intersects $\Gamma_L$ with a certain unitary group:
\smallskip
The lattice $L$ admits a complex structure. We follow Kondo's description
[Ko2].
For this one uses the isomorphism $\gz^4\to\gz^4$ defined by the matrix
$$\pmatrix{-1&0&2&0\cr 0&1&0&2\cr -1&0&1&0\cr 0&-1&0&-1\cr}.$$
This is an isometry
$$\rho_1:U\oplus U(2)\Isom U\oplus U(2)$$
We also use the isometry
$$\rho_0:D_4\Isom D_4,\quad (x_1,x_2,x_3,x_4)\loma (x_2,-x_1,x_4,-x_3).$$
The direct sum $\rho:=\rho_1\oplus\rho_0\oplus\rho_0$ defines an isometry
$$\rho:L\Isom L.$$
Since $\rho^2=-\id$
we get a structure as $\gz[\imag]$ module on $L$ by defining
$$\imag a:=\rho(a).$$
One has $L\cong \gz[i]^6$.
We denote by $(a,b)=q(a+b)-q(a)-q(b)$ the real bilinear
form on $L$. Because of  $(a,\rho(a))=0$ the pairing
$$2\spitz{a,b}:=(a,b)-\imag(a,\imag b)\qquad
(\Longrightarrow 2\spitz{a,a}=(a,a)).$$
is a hermitian form of signature $(1,5)$.
Notice that we take hermitian forms to be
$\cz$-linear in the second variable.
Let $\U(L,\spitz{\cdot,\cdot})$ denote the unitary group. We set
$$\Gamma:=\O'(L)\cap\U(L,\spitz{\cdot,\cdot}),\quad
\Gamma[1-\imag]:= \Gamma_L\cap\U(L,\spitz{\cdot,\cdot}).$$
The natural homomorphism
$$\Gamma/\Gamma[1-\imag]\lo \O'(L)/\Gamma_L\lo\O(\fz_2^6)\cong S_8$$
is still an isomorphism.
\smallskip
We denote by $\calB$ the associated complex five dimensional ball.
It can be considered as the set of positive definit
one-dimensional complex subspaces
of $L\otimes_\gz\rz$. We can intersect the Heegner divisor $\calH_\alpha$
with the ball:
$$\calB_\alpha:=\calH_\alpha\cap\calB.$$
The group $\Gamma$ acts on $\calB$  and one can consider the graded algebras
of modular forms $A(\Gamma[1-\imag])$.
There is a natural restriction homomorphism $A(\Gamma_L)\to A(\Gamma[1-\imag])$.
Together with \Bal\ we get a homomorphism $\calR_3/\calI_3\to A(\Gamma[1-\imag])$.
\proclaim
{Theorem}
{The homomorphism
$$\calR_3/\calI_3\to A(\Gamma[1-\imag])$$
is an isomorphism. Hence this ring of modular
forms is generated by\/ $14$ modular forms
with\/ $14$ defining quadratic relations. The dimension formula is given by \HF.
\smallni
{\bf Corollary.} There is an induced biholomorphic map
$$\bar X(8)\lo \overline{\calB/\Gamma[1-\imag]} ,$$
where the right hand side denotes the Baily-Borel compactification.
}
Rmf%
\finishproclaim
This result can be found in Kondo's paper [Ko2]. But the proof is not convincing.
Kondo constructs by means of the theory of $K3$-surfaces a period map
$X(8)\lo {\calB/\Gamma[1-\imag]}$. Then he wants to use an extension theorem
of Borel to extend it to the compactifications. This extension theorem of
Borel states
essentially that a holomorphic map of a punctured disk into an arithmetic quotient
$D/\Gamma$ extends to a holomorphic map of the full disk to the
Satake compactification.
The theorem of Borel requires that $\Gamma$ is torsion free. Otherwise it is false
as for example $\SL(2,\gz)$ shows. Also $\Gamma[1-\imag]$
can be taken as counter example.
\smallskip
The group
$\Gamma[1-\imag]$ is not torsion free.
If $r$ is an element with $2\spitz{r,r}=(r,r)=-2$ then the complex reflection
$$x\loma x+2\spitz{x,r}x=x+(x,r)x-\imag(x,\imag r)x$$
is contained in $\Gamma[1-\imag]$.
These reflections are related to the divisors $\calB_\alpha=\calH_\alpha\cap\calB$.
By definition this divisor consists of all $z\in\calB$ such that
$$(\delta,z)=(\delta,\imag z)=0\quad\hbox{for all}\quad \delta\in L^*,\quad
(\delta,\delta)=-1,
\quad\hbox{over}\quad \alpha.$$
Using the hermitian form this is equivalent to $\spitz{\delta,z}=0$.
This we can rewrite as
$$\spitz{r,z}=0,\qquad r:=(1-\imag)\delta.$$
Obviously $r\in L$ and $\spitz{r,r}=-1$. This shows that $\calB_\alpha$ is the fixed point
set of a reflection inside $\Gamma[1-\imag]$.
\proclaim
{Remark}
{The zero divisor of the restriction $F_S=f_S\vert\calB$ of the star modular form
is contained in the ramification divisor of
the projection $\calB\to\calB/\Gamma[1-\imag]$.}
Rsr%
\finishproclaim
Because of the ramification the
star modular forms $F_S$, restricted to the ball,
vanish at $\calB_\alpha$, $\alpha\in S$ of order two.
The ramification makes another correction in [Ko2]
necessary. It is possible to take a holomorphic square root $G_S$ of $F_S$.
These forms are
modular forms with characters. Hence for two stars $S_1,S_2$ the quotient
$G_{S_1}/G_{S_2}$
is not invariant under $\Gamma[1-\imag]$.
The stars and the tableaux are in 1-1 correspondence [Ko2].
Denote by $\tau_1$, $\tau_2$
the tableaux which correspond to $S_1$, $S_2$ and by $\mu_{\tau_1}$, $\mu_{\tau_2}$
the corresponding Specht polynomials. Kondo claims
$G_{S_1}/G_{S_2}=\mu_{\tau_1}/\mu_{\tau_2}$
(s.~the two lines before theorem 7.6 in [Ko2]). The correct equation is
$$F_{S_1}/F_{S_2}=\mu_{\tau_1}/\mu_{\tau_2}.$$
The point is that the pole- and zero
orders of $F_{S_1}/F_{S_2}$ considered as functions
on the {\it quotient\/} $\calB/\Gamma[1-\imag]$ are only one and not two.
\smallni
{\it Proof of theorem \Rmf.\/} From the description of the zeros of the $F_S$
it follows that they have no common zero on the
Baily-Borel compactified ball quotient. Hence one
obtains a finite map
$$\overline{\calB/\Gamma[1-\imag]}\lo\proj(\calR_3/\calI_3).$$
Since $\calR_3/\calI_3$ is normal,  one has only to show that this map is
generically injective.
This follows from Kondo's comparison from $X(8)$ and $\calB/\Gamma[1-\imag]$
by means of K3-surfaces.
It can also be proven directly using the methods of [FS]. In this way one can
avoid the use of K3-surfaces.
(It depends on the situation and also on a question of taste whether one
wants to apply
the theory of moduli to the theory of modular forms or conversely.)
\neupara{Siegel modular forms of genus three}%
Let
$$\calH_g=\set{Z\in \cz^{(g,g)};\quad Z={^tZ},\ \Im Z>0}$$
be the Siegel half plane of genus $g$ and
$$\Gamma_g[q]=\kernel\bigl(\Sp(g,\gz)\lo \Sp(g,\gz/q\gz)\bigr)$$
the principal congruence subgroup of level $q$ in the Siegel modular group
$\Gamma_g:=\Gamma_g[1]$. It acts on $\calH_g$ by the usual formula
$MZ=(AZ+B)(CZ+D)^{-1}$. The quotient $\calH_g/\Gamma_g[q]$ is, by a well-known
theorem of Baily, a quasi projective variety. We denote by
$\overline{\calH_g/\Gamma_g[q]}$ its Satake compactification.
This is related to the algebra of modular forms
$$A(\Gamma_g[q])=\bigoplus_{r=0}^\infty[\Gamma_g[q],r],$$
where $[\Gamma_g[q],r]$ denotes the vector space of modular forms of weight
$r$. These are holomorphic functions $f:\calH_g\to\cz$ with the transformation property
$f(MZ)=\det(CZ+D)^rf(Z)$ for all $M\in \Gamma_g[q]$ (and a regularity condition at
the cusps in the case $g=1$).
Then one has
$$\overline{\calH_g/\Gamma_g[q]}\cong\proj(A(\Gamma_g[q])).$$
We denote by $\calJ_g\subset\calH_g$ the set of all matrices
which are period matrices of hyperelliptic (non-singular projective) curves.
Then
$$\calI_g[q]=\calJ_g/\Gamma_g[q]$$
is a quasi projective subvariety.
The variety $\calI_g[1]$ is the modular variety of hyperelliptic curves.
The closure
$$\overline{\calI_g[q]}\subset\overline{\calH_g/\Gamma_g[q]}$$
is one of its compactifications we
have to consider.
\smallskip
For $g\geq 3$ and $q\geq 2$, $\overline{\calI_g[q]}$ is reducible. The number of reducible components  of  $\overline{\calI_g[2]}$ is  $$2^{2g^2 +g}\prod_{k=1}^g(1-2^{-2k})/(2g+2)!,$$ cf.
[Ts] and  it is the same for all even $q$.
To describe the irreducible
components we have to recall the thetanullwerte which for arbitrary $g$ are defined
by
$$\vartheta[m]=\sum_{n\in\gz^g}e^{\pii\bigl(Z[n+m'/2]+(n+m'/2)'m''\bigr)},\qquad
m=\pmatrix{m'\cr m''}.$$
Here $m',m''$ are columns in $\gz^g$. Up to the sign they depend only on $m',m''$ mod 2.
They vanish identically if and only if $^tm'm''$ is odd. Hence in the case $g=3$,  there are
essentially 36 such thetas  and 36 irreducible components of $\overline{\calI_3[2]}$ .
By classical results a matrix $Z\in\calH_3$ is in the closure of $\calJ_3$ if and only if the
product of all 36 thetas vanishes. From this one deduces that
$\calI_3[2]$ has 36 irreducible components. Each of them corresponds to
a so-called ``even
characteristic''
$$m=\pmatrix{m'\cr m''}\in\fz_2^6,\quad ^tm'm''=0,$$
and is defined as
$$\calI_3[2]^m:=\set{Z\in \calI_3[2];\quad \vartheta[m](Z)=0}.$$
Here we remark that $\vartheta[m]^4$ is a modular form of weight 2 on $\Gamma_2[2]$.
Hence the zero locus of $\vartheta[m]$
inside $\calI_3[2]$ is well defined. We have to consider the closure
$$\overline{\calI_3[2]}^m\subset \overline{\calH_3/\Gamma_g[2]}.$$
The different components are permuted under the full modular group.
Hence it doesn't matter, which characteristic we use. We can take for example the
zero characteristic $0$.
The Galois group of the covering
$$\overline{\calI_3[2]}^m\lo \overline{\calI_3[1]}$$
is the stabilizer of $m$ in $\Sp(2,\gz/2\gz)$
with respect to the usual action of this group on characteristics:
$$M\{m\}:={^tM}^{-1}m+\pmatrix{(A\,{^tB})_0\cr (C\,{^tD})_0}\qquad (\mod\; 2).$$
Here we use the notation $A_0$ for the column vector built from the diagonal of
a matrix $A$.
In case of the zero characteristic the stabilizer is
the subgroup of $\Sp(2,\fz_2)$, which fixes the
quadratic form
$$\fz_2^6\lo\fz_2,\qquad m\loma {^tm'm''}.$$
Hence it is our orthogonal group $\O(\fz_2^6)$.
This is also the image of the theta group $\Gamma_{3,\vartheta}$, where
$$\Gamma_{g,\vartheta}:=\set{M\in \Gamma_g;\quad
(A\,{^tB})_0\equiv (C\,{^tD})_0\equiv 0\ \mod 2}.$$
Hence we have
$$\O(\fz_2^6)=\Gamma_{3,\vartheta}/\Gamma_3[2].$$
Next we have to make use of the theory of Thomae, who computed the
theta values of hyperelliptic curves. We use (slightly modified) the beautiful
approach of Mumford [Mu]. Let be
$$
B=\{1,\dots,8\},\quad
U=\{1,\dots,4\}.$$
Recall that we consider on $\fz_2^6$ the quadratic form
$$q:\fz_2^6\lo\fz_2,\quad q(m)={^tm'm''}$$
and the associated pairing
$$(m,n)=q(m+n)+q(m)+q(n)={^tm'n''}+{^tm'' n'}.$$
Their multiplicative form is
$$e(m)=(-1)^{^tm'm''},\quad
e(m,n)=(-1)^{^tm'n''+^tm'' n'}.$$
Recall that a characteristic $m$ is called even if $e(m)=1$.
Hence the even charactersitics are just the isotropic vectors of the quadratic
space $(\fz_2^6,q)$.
There are 36 even $m$.
Following  Mumford [Mu], we construct a map
$$T\loma m(T),$$
which attaches to an subset $T\subset B$ of even order an
element from $\fz_2^6$.
The basic properties of this assignment are
\vskip2mm
\item{a)} $m(T_1)=m(T_2)$  if and only if $T_1=T_2$  or
$T_1$ is the complement of $T_2$. Especially every element of $\fz_2^6$
is in the image of the map $m$.
\vskip1mm
\item{b)} With the notation
$$T_1\circ T_2=T_1\cup T_2- T_1\cap T_2$$
one has
$$m(T_1\circ T_2)=m(T_1)+m(T_2) +m (\emptyset).$$

\vskip1mm
\item{c)} The characteristic $m(T)$ is even if and only if the number
of elements of $T\circ U$ is divisible by 4 (hence it is $0$, $4$ or $8$).
\smallni
 \item{d)}\hfil
$\displaystyle e(m(T_1))e(m(T_2))e(m(T_1\circ T_2))=(-1)^{\# T_1\cap T_2}  $.\hfil
\medni
As consequence  we get

$$e(m(T_1),m(T_2))= e(m(T_1))e(m(T_2)) e(m(T_1)+(m(T_2))=$$
$$e(m(T_1))e(m(T_2))e(m(T_1\circ T_2)+m(\emptyset))=$$
$$e(m(T_1))e(m(T_2))e(m(T_1\circ T_2))e(m(\emptyset))e(m(T_1\circ T_2),m(\emptyset))=$$
$$(-1)^{\# T_1\cap T_2}e(m(\emptyset)) e(m(T_1\circ T_2), m(\emptyset))$$

We observe that if we set $m(\emptyset)$ equal to the zero characteristics ,
then we have Mumford's
procedure, but this will not be our case.
To be concrete we take the following explicit realization.
\bigni
\centerline{\vbox{
\halign{$\hfil\ #\ \hfil$&$\hfil\ #\ \hfil$&$\hfil\ #\ \hfil$
&$\hfil\ #\ \hfil$&$\hfil\ #\ \hfil$
&$\hfil\ #\ \hfil$&$\hfil\ #\ \hfil$&\quad$\hfil#\hfil$\cr
\{1,2\}&\{1,3\}&\{1,4\}&\{1,5\}&\{1,6\}&\{1,7\}&\{1,8\}&\emptyset\cr
\noalign{\vskip3pt}
1&1&0&0&1&0&1&0\cr 1&1&0&1&0&0&0&1\cr 1&0&1&1&0&0&1&0\cr
\noalign{\vskip1pt}
1&0&1&0&0&1&0&1\cr 0&1&1&0&0&1&1&0\cr0&1&1&0&1&0&0&1\cr}}}
\medni

The first row contains subsets of $B$, the columns below are the corresponding
characteristics. It is easy to verify that this assignment uniquely extends
to a map with the properties a)-d). We also notice that $m(U)$ is the sum of the
first three,
$$^tm(U)=(0,0,0,0,0,0).$$
The group $\Sp(3,\gz/2\gz)$ acts on the set of characteristics, we recall from [Ig1]
that  $e(m)$ and $e(m,n)e(m,p) e(n, p) $ are invariant.
The form $e(m,n)$ results to be  $\O(\fz_2^6)$ invariant.
\smallskip
The group $S_8$ acts on $B=\{1,\dots,8\}$ and hence on the pairs $\{T,B-T\}$
of subsets of even order and
on $U$,
so we get an induced action of $S_8$
on $\fz_2^6$.
We observe that  $m(U)=0$  and $ \# (T \circ U)$ are invariant for the action of $S_8$, in fact
it generates   the subgroup $\O(\fz_2^6)$ of $\GL(\fz_2^6)$.
Hence we have constructed now an explicit isomorphism
$$S_8\Isom \O(\fz_2^6).$$
This isomorphism will be used in the rest of this paper.
We associate to each subset $T\subset B$ of even order a monomial
$D(T)$
in the ring
$\cz[X_1,\dots,X_8]$.
We have to distinguish two cases:
\smallni
1) $T\circ U\ne 4$.  In this case we attach $0$.
\smallni
2) $T\circ U=4$. In this case we attach the following monomial of degree
12
$$e(m(T), m(\emptyset))(-1)^{\# T\cap U}\Delta\prod_{i\in T\circ U\atop j\notin T\circ U}
(X_i-X_j)^{-1},\qquad
\hbox{where}\quad
\Delta=\prod_{i<j}(X_i-X_j).$$
Because of $D(T)=D(B-T)$ the monomial
$$D(m):=D(m(T))\quad\hbox{for}\quad m=m(T)$$
is well defined.
\smallskip
From Thomae's  computation of the theta values of hyperelliptic
Riemann surfaces follows (s.~[Mu]):
\proclaim
{Theorem (Thomae, 1870)}
{Let $\cz[\dots\vartheta[m]^4\dots]$ be the ring generated by the fourth powers
of the 36 even thetanullwerte. Then $D$ defines a ring homomorphism
$$D:\cz[\dots\vartheta[m]^4\dots]\lo \cz[\dots X_i-X_j\dots].$$
This homomorphism is equivariant with the standard action of $\O(\fz_2^6)$
on the thetanullwerte and the action of $S_8$ on the variables $X_i$.
}
ST%
\finishproclaim
The only element of the generators which lies in the kernel is
$\vartheta[0]^4$. But this does of course not mean that the kernel is generated
by this element.
\smallskip
We investigate the image. Let $m$ be  the characteristic coming from the empty set, i.e
$^tm=(0,1,0,1,0,1)$, then
 we have
$$\eqalign{
D(\vartheta[m]^4)&=(X_1-X_2)(X_1-X_3)(X_1-X_4)(X_2-X_3)(X_2-X_4)(X_3-X_4)\cr
&(X_5-X_6)(X_5-X_7)(X_5-X_8)(X_6-X_7)(X_6-X_8)(X_7-X_8).\cr}
$$
This is the product of the three Specht polynomials with respect to the tableaux
$$\pmatrix{1&3&5&7\cr 2&4&6&8}\quad
\pmatrix{1&2&5&6\cr 3&4&7&8}\quad
\pmatrix{1&2&5&6\cr 4&3&8&7}
$$
This  homomorphism can be
given explicitly as
\proclaim
{Proposition}
{There is a homomorphism
$$\xymatrix{
  D:\cz[\dots\vartheta[m]^4\dots]\ar[rd]\ar[r]&\cz[X_1,\dots,X_8] \\
&\cz[Y_1,\dots,Y_{14}]/\calJ\ar@^{(->}[u]\ar@{->}[u]}
$$
Its kernel contains $\vartheta[0]^4$. An explicit description, using the notation
$W_{ij}=X_i-X_j$, is
given in the following table.}
VExH%
\finishproclaim
\halign{\quad$#$\hfil&\quad$\hfil#\;$&$#$\hfil\cr
(0,0,0,0,0,0)&&
0\cr
(0,0,0,0,0,1)&
-& W_{13}W_{16}W_{17}W_{24}W_{25}W_{28}W_{36}W_{37}W_{45}W_{48}W_{58}W_{67}\cr
(0,0,0,0,1,0)&
 -&W_{12}W_{13}W_{15}W_{23}W_{25}W_{35}W_{46}W_{47}W_{48}W_{67}W_{68}W_{78}\cr
(0,0,0,0,1,1)&
&W_{13}W_{14}W_{18}W_{25}W_{26}W_{27}W_{34}W_{38}W_{48}W_{56}W_{57}W_{67}\cr
(0,0,0,1,0,0)&
-&W_{14}W_{15}W_{17}W_{23}W_{26}W_{28}W_{36}W_{38}W_{45}W_{47}W_{57}W_{68}\cr
(0,0,0,1,0,1)&
&W_{12}W_{17}W_{18}W_{27}W_{28}W_{34}W_{35}W_{36}W_{45}W_{46}W_{56}W_{78}\cr
(0,0,0,1,1,0)&
&W_{15}W_{16}W_{18}W_{23}W_{24}W_{27}W_{34}W_{37}W_{47}W_{56}W_{58}W_{68}\cr
(0,0,0,1,1,1)&
-&W_{12}W_{14}W_{16}W_{24}W_{26}W_{35}W_{37}W_{38}W_{46}W_{57}W_{58}W_{78}\cr
(0,0,1,0,0,0)&
&W_{13}W_{14}W_{16}W_{25}W_{27}W_{28}W_{34}W_{36}W_{46}W_{57}W_{58}W_{78}\cr
(0,0,1,0,1,0)&
-&W_{13}W_{17}W_{18}W_{24}W_{25}W_{26}W_{37}W_{38}W_{45}W_{46}W_{56}W_{78}\cr
(0,0,1,1,0,0)&
-&W_{12}W_{14}W_{18}W_{24}W_{28}W_{35}W_{36}W_{37}W_{48}W_{56}W_{57}W_{67}\cr
(0,0,1,1,1,0)&
&W_{12}W_{16}W_{17}W_{26}W_{27}W_{34}W_{35}W_{38}W_{45}W_{48}W_{58}W_{67}\cr
(0,1,0,0,0,0)&
&W_{12}W_{15}W_{16}W_{25}W_{26}W_{34}W_{37}W_{38}W_{47}W_{48}W_{56}W_{78}\cr
(0,1,0,0,0,1)&
-& W_{14}W_{16}W_{18}W_{23}W_{25}W_{27}W_{35}W_{37}W_{46}W_{48}W_{57}W_{68}\cr
(0,1,0,1,0,0)&
-&W_{13}W_{15}W_{18}W_{24}W_{26}W_{27}W_{35}W_{38}W_{46}W_{47}W_{58}W_{67}\cr
(0,1,0,1,0,1)&
&W_{12}W_{13}W_{14}W_{23}W_{24}W_{34}W_{56}W_{57}W_{58}W_{67}W_{68}W_{78}\cr
(0,1,1,0,0,0)&
&W_{16}W_{17}W_{18}W_{23}W_{24}W_{25}W_{34}W_{35}W_{45}W_{67}W_{68}W_{78}\cr
(0,1,1,0,1,1)&
&W_{12}W_{15}W_{18}W_{25}W_{28}W_{34}W_{36}W_{37}W_{46}W_{47}W_{58}W_{67}\cr
(0,1,1,1,0,0)&
-&W_{12}W_{13}W_{17}W_{23}W_{27}W_{37}W_{45}W_{46}W_{48}W_{56}W_{58}W_{68}\cr
(0,1,1,1,1,1)&
-&W_{13}W_{15}W_{16}W_{24}W_{27}W_{28}W_{35}W_{36}W_{47}W_{48}W_{56}W_{78}\cr
(1,0,0,0,0,0)&
&W_{12}W_{14}W_{17}W_{24}W_{27}W_{35}W_{36}W_{38}W_{47}W_{56}W_{58}W_{68}\cr
(1,0,0,0,0,1)&
-& W_{15}W_{17}W_{18}W_{23}W_{24}W_{26}W_{34}W_{36}W_{46}W_{57}W_{58}W_{78}\cr
(1,0,0,0,1,0)&
-&W_{12}W_{16}W_{18}W_{26}W_{28}W_{34}W_{35}W_{37}W_{45}W_{47}W_{57}W_{68}\cr
(1,0,0,0,1,1)&
&W_{14}W_{15}W_{16}W_{23}W_{27}W_{28}W_{37}W_{38}W_{45}W_{46}W_{56}W_{78}\cr
(1,0,1,0,0,0)&
&W_{14}W_{15}W_{18}W_{23}W_{26}W_{27}W_{36}W_{37}W_{45}W_{48}W_{58}W_{67}\cr
(1,0,1,0,1,0)&
 -&W_{15}W_{16}W_{17}W_{23}W_{24}W_{28}W_{34}W_{38}W_{48}W_{56}W_{57}W_{67}\cr
(1,0,1,1,0,1)&
&W_{13}W_{16}W_{18}W_{24}W_{25}W_{27}W_{36}W_{38}W_{45}W_{47}W_{57}W_{68}\cr
(1,0,1,1,1,1)&
-&W_{13}W_{14}W_{17}W_{25}W_{26}W_{28}W_{34}W_{37}W_{47}W_{56}W_{58}W_{68}\cr
(1,1,0,0,0,0)&
&W_{12}W_{13}W_{18}W_{23}W_{28}W_{38}W_{45}W_{46}W_{47}W_{56}W_{57}W_{67}\cr
(1,1,0,0,0,1)&
-& W_{13}W_{14}W_{15}W_{26}W_{27}W_{28}W_{34}W_{35}W_{45}W_{67}W_{68}W_{78}\cr
(1,1,0,1,1,0)&
&W_{14}W_{17}W_{18}W_{23}W_{25}W_{26}W_{35}W_{36}W_{47}W_{48}W_{56}W_{78}\cr
(1,1,0,1,1,1)&
-&W_{12}W_{15}W_{17}W_{25}W_{27}W_{34}W_{36}W_{38}W_{46}W_{48}W_{57}W_{68}\cr
(1,1,1,0,0,0)&
&W_{13}W_{15}W_{17}W_{24}W_{26}W_{28}W_{35}W_{37}W_{46}W_{48}W_{57}W_{68}\cr
(1,1,1,0,1,1)&
&W_{12}W_{13}W_{16}W_{23}W_{26}W_{36}W_{45}W_{47}W_{48}W_{57}W_{58}W_{78}\cr
(1,1,1,1,0,1)&
&W_{14}W_{16}W_{17}W_{23}W_{25}W_{28}W_{35}W_{38}W_{46}W_{47}W_{58}W_{67}\cr
(1,1,1,1,1,0)&
&W_{12}W_{14}W_{15}W_{24}W_{25}W_{36}W_{37}W_{38}W_{45}W_{67}W_{68}W_{78}\cr
}
\proclaim
{Remark}
{The variable $W_{ij}$ appears in the monomial $D(\vartheta[m]^4)$, $m\ne 0$,
if and only if $m$ is orthogonal to the (odd) characteristic
$m(T_{ij})$, where $T_{ij}$ is defined by
$$T_{ij}\circ U=\{i, j\}.$$
}
RSf%
\finishproclaim
It is easy to verify this for one $m$ by direct computation. The general case follows
then because $\O(\fz_2^6)$ acts transitively on the set of non-zero even
characteristics.\qed\smallskip
Finally we mention an observation of Igusa [Ig2]:
\proclaim
{Remark}
{ $D$ extends to  a ring homomorphism
$$D: A(\Gamma_3[2])\lo\cz[Y_1,\dots,Y_{14}]/\calJ.$$
An element $f$ belongs to the kernel if and only if it vanishes on $\calI_3[2]^0$.}
SU%
\finishproclaim
Recall that on the right hand side we have the ring of Specht polynomials from section 1,
which also can be identified with
the ring of ball modular forms $\calR_3/\calI_3$ from section 2.
\neupara{A structure theorem}%
Runge [Ru2] has described generators for the ring of Siegel modular forms of genus
three and level two with respect to the trivial character. He gave explicitly
15 generators in weight 2 and 15 generators in weight 4.
He expressed them as explicit polynomials in the 8 theta
constants of second kind
$$f_a(Z)=\sum_{g\in\gz^3}e^{2\pii Z[g+a/2]},\qquad a\in (\gz/2\gz)^3.$$
We recall the relation to the $\vartheta[m]$:
$$
\vartheta^2\Bigl[{m'\atop m''}\Bigr]=\sum_{a\in\fz_2^3}
(-1)^{^t m' m''}f_{m'+a}f_a,\quad
f_af_b={1\over 8}\sum_{c\in\fz_2^3}(-1)^{^tac}
\vartheta^2\Bigl[{a+b\atop c}\Bigr].$$
The classical Schottky relation $R$ can be written as a polynomial in degree
$16$ of the $f_a$. Runge proved the basic result that
this is a defining relation in the ring
$\cz[f_a]$
and moreover that this ring is normal. This implies that the subring
$\cz[f_af_b]$ of forms of integral weight is the full ring of modular
forms of integral weight and trivial character with respect to Runge's
group $\Gamma_3^*[2,4]$. This is the unique subgroup of Igusa's group
$\Gamma_3[2,4]$ which dos not contain the negative unit matrix.
Recall  that the Igusa group $\Gamma_g[q,2q]$ is the subgroup
of all $M\in\Gamma_g[q]$ such that $A\,{^tB}/q$ and $C\,{^tD}/q$ have even diagonal.
Hence the theta group is $\Gamma_{g,\vartheta}=\Gamma_g[1,2]$.
\smallskip
We recall briefly the way how Runge describes the action of the modular group
on the $f_a$. He introduces the following group $H_3$ of $8\times 8$-matrices
(indexed by the 8 elements from $\fz_2^3$).
The first generator is
$$\tilde S=\Bigl({1+\imag\over2}\Big)^3\bigl((-1)^{^tab}\bigr).$$
To every integral diagonal $3\times 3$-matrices one attaches the diagonal matrix
$\tilde T_S$ with diagonal entries $\imag^{S[a]}$.
The group $H_3$ is generated by $\tilde I$ and by the $\tilde T_S$.
The negative unit matrix is contained in $H_3$. If we map $I$ to $\tilde I$ and
$T_S$ to $\tilde T_S$ we obtain a homomorphism
$$\Gamma_3\lo H_3/\pm.$$
the kernel of this homomorphism is $\Gamma^*[2,4]$.
Moreover there is a surjective homomorphism
$$H_3\lo\Sp(3,\fz_2),$$
which sends $\tilde I$ and $\tilde T_S$ to the cosets of $I$ and $T_S$.
The kernel of this homomorphism is
$$N_3:=\langle \imag E, \tilde T_{2S}, \tilde I^{-1}\tilde T_{2S} I,\quad
(S\ \hbox{integral})\rangle,$$
as correctly stated in Runge's first paper [Ru1].
In the second part of his paper [Ru2] Runge claims that $N_3$ agrees with
$$N_3':=\langle\tilde T_{2S}, \tilde I^{-1}\tilde T_{2S} I,\quad
(S\ \hbox{integral})\rangle.$$
But this is not the case, since $N_3'\subset N_3$ is a subgroup of index
two which doesn't contain $\imag E$. It contains only $-E$.
\smallskip
Now we modify Runge's results on modular forms of level two introducing the
theta character $v_\vartheta$. This is the character of the modular form
$\vartheta[0]^2$, which lives on the theta group $\Gamma_{3,\vartheta}=
\Gamma_3[1,2]$. This group contains the principal congruence group
of level two
$\Gamma_3[2]$. We consider the ring
$$A(\Gamma_3[2]):=\bigoplus_{r\in\gz}[\Gamma_3[2],r,v_\vartheta^r].$$
We can get now a formal description of this ring:
\proclaim
{Remark}
{let $F_a$ be a formal variable for $a\in\fz_2^3$. We consider the
polynomial ring\/  $\cz[F_a]$ and the polynomial $R$ of degree $16$ which
describes the Schottky relation. Since $R$ is invariant under $H_3$ the group
$N_3'$ acts on $\cz[F_a]/(R)$. There is a natural isomorphism
$$\bigl(\cz[F_a]/(R)\bigr)^{N_3'}=
\bigl(\cz[F_a]^{N_3'}\bigr)/(R)
\Isom A(\Gamma_3[2]),$$
which replaces $F_a$ by $f_a$.}
RMf%
\finishproclaim
From Runge's results  one can derive:
\proclaim
{Proposition}
{The ring $A(\Gamma_3[2])$ is generated by
\smallni
$1$ form of weight\/ $1$ (namely $\vartheta[0]^2$),
which belongs to the theta group,\hfill\break
$14$ forms of weight\/ $2$,\hfill\break
$14$ forms of weight\/ $3$,\hfill\break
$1$ form of weight\/  $4$, which belongs to the theta group.
\smallni
The Hilbert function of the ring is  $(1-z^8)$ times
$${1-3z+13z^2-17z^3+44z^4-17z^{5}+13z^{6}-3z^{7}+z^{8}\over
(1-z^2)^4(1-z)^4}\;,$$
which is
$$\eqalign{
&1+z+15z^2+29z^3+135z^4+310z^5+870z^6+1830z^7+\cr&3992z^8+7534z^9+14142z^{10}+\cdots\cr}
$$
}
RA%
\finishproclaim
\neupara{The ring of hyperelliptic modular forms of genus 3}%
The factor group $\Gamma_3/\Gamma_3[2]$ is isomorphic to $\Sp(3,\gz/2\gz)$.
This group contains naturally the orthogonal group $\O(\fz_2^6)$
(quadratic form $m_1m_4+m_2m_5+m_3m_6$).
Recall that the inverse image of $\O(\fz_2^6)$ in the Siegel modular group is
the theta group $\Gamma_{3,\vartheta}$. The group $\O(\fz_2^6)$ has a
simple subgroup of index two,
which we denote by $\SO(\fz_2^6)$. Its inverse image is a subgroup of index two
$\Gamma_{3,\vartheta}'\subset \Gamma_{3,\vartheta}$.

It is important to find a natural system of generators of the  space
$[\Gamma_3[2],2]$. Actually we will construct 15 forms $\Theta_1,\dots,\Theta_{15}$,
which generate this space and which have the
property that they are permuted under the action of the  group $\Gamma_{3,\vartheta}'$.
For their construction we use the 36 theta constants of first kind
$\vartheta[m]$.
\smallskip
We have thirty maximal totally isotropic spaces  contained  in $\fz_2^6$.
For such a  subspace $A$ we set
$$\Theta=\sum_{i=1}^7\vartheta[m_i]^4,$$
where $m_1,\dots,m_7$ are the non-zero elements of $A$.
The orthogonal group permutes  the totally isotropic spaces,
the simple group $\SO(\fz_2^6)$  describes two orbits.
(Two spaces  are in the same orbit if and only if their intersection has even dimension.)
We give the
list of the subspaces in one of the two orbits.
We  give  their entries  different from $0$ ordered
lexicographically, i.e. $m$ is replaced by the digit $2^5m_1+2^4m_2+2^3m_3+2^2m_4+2 m_5 +m_6$.
The fifteen orbits are:
\medni
\vbox{
\halign{\qquad$#$\hfil&$=#$\hfil&\quad$#$\hfil&$=#$\hfil\cr
A_1&[6,8,14,48,54,56,62]&A_2&[6,10,12,49,55,59,61]\cr
A_3&[7,24,31,40,47,48,55]&A_4&[7,27,28,,42,45,49,54]\cr
A_5&[4,8,12,16,20,24,28]&A_6&[4,10,14,17,21,27,31]\cr
A_7&[5,16,21,40,45,56,61]&A_8&[5,17,20,42,47,59,62]\cr
A_9&[2,8,10,32,34,40,42]&A_{10}&[2,12,14,33,35,45,47]\cr
A_{11}&[3,24,27,32,35,56,59]&A_{12}&[3,28,31,33,34,61,62]\cr
A_{13}&[1,16,17,32,33,48,49]&A_{14}&[1,20,21,34,35,54,55]\cr
A_{15}&[1, 2,3,4,5,6,7]\cr
}}
\medni
For  each $A_i$ we denote by $\Theta_i$ the corresponding sum of
fourth powers of the thetanullwerte. For example with these notations we have

$$\leqalignno{
\Theta_{1}=
&\phantom{+\;}\vartheta\left[\matrix{0&1\cr0&1\cr0&0\cr}\right]^4
+\vartheta\left[\matrix{0&0\cr0&0\cr1&0\cr}\right]^4
+\vartheta\left[\matrix{0&1\cr0&1\cr1&0\cr}\right]^4\cr&
+\vartheta\left[\matrix{1&0\cr1&0\cr0&0\cr}\right]^4
+\vartheta\left[\matrix{1&1\cr1&1\cr0&0\cr}\right]^4
+\vartheta\left[\matrix{1&0\cr1&0\cr1&0\cr}\right]^4
+\vartheta\left[\matrix{1&1\cr1&1\cr1&0\cr}\right]^4.
\cr}$$

\proclaim
{Proposition}
{ The forms $\Theta_1,\dots,\Theta_{15}$ give a basis of the space
$[\Gamma_3[2],2]$.
As a consequence of the quartic Riemann relations the form
$$\Theta_1+\cdots+\Theta_{15}$$
vanishes on the hyperelliptic component defined by $\vartheta[0]=0$.
}
Fun%
\finishproclaim
We are interested in the restriction of the ring $A(\Gamma_3[2])$
to one of the hyperelliptic components. We take the component defined by the zero locus
of $\vartheta[0]$. It turns out that not only $\vartheta[0]^2$ but also the
other $\Gamma_{3,\vartheta}$ invariant form in \RA\ of weight 4 vanishes on this component.
We give them explicitly as polynomials in the $F_a$. We order them
lexicographically, i.e. $a$ is replaced by the digit $4a_1+2a_2+a_3$.
\medni
{\ninepoint
$P: =
F_0^2+F_1^2+F_2^2+F_3^2+F_4^2+F_5^2+F_6^2+F_7^2$,
\smallni
$Q:=
F_2^8+4F_4^4F_2^2F_3^2+F_1^8-2F_7^4F_0^4+4F_6^2F_7^2F_2^4-2F_7^4F_1^4+
24F_6^2F_7^2F_2^2F_3^2-8F_0^2F_2^2F_5^2F_6^2-4F_5^6F_4^2-2F_3^4F_4^4+
F_5^8+4F_1^4F_2^2F_3^2-2F_1^4F_2^4-2F_2^4F_0^4-2F_0^4F_3^4-10F_1^4F_0^4-
4F_1^6F_0^2+4F_0^4F_3^2F_2^2-2F_6^4F_0^4-2F_3^4F_1^4+4F_3^4F_1^2F_0^2+
4F_5^4F_2^2F_3^2-2F_5^4F_1^4+4F_5^4F_0^2F_1^2+4F_6^4F_3^2F_2^2-2F_7^4F_2^4-
2F_7^4F_3^4-2F_6^4F_1^4-2F_6^4F_3^4-2F_6^4F_2^4-2F_5^4F_0^4-2F_5^4F_2^4-
2F_5^4F_3^4+4F_7^4F_2^2F_3^2+4F_0^2F_1^2F_2^4+4F_6^4F_1^2F_0^2+
4F_7^4F_0^2F_1^2-128F_1F_0F_3F_2F_5F_4F_7F_6-10F_4^4F_5^4-4F_6^6F_7^2-
10F_6^4F_7^4-4F_6^2F_7^6+4F_4^2F_6^4F_5^2+4F_4^4F_7^2F_6^2+4F_5^4F_6^2F_7^2-
2F_4^4F_6^4-2F_4^4F_7^4-2F_5^4F_6^4+F_6^8+F_4^8+F_7^8-2F_5^4F_7^4-4F_4^6F_5^2+
F_0^8+F_3^8+4F_0^2F_1^2F_4^4+4F_4^2F_5^2F_7^4-2F_2^4F_4^4-4F_0^6F_1^2-
4F_2^6F_3^2-10F_2^4F_3^4-4F_2^2F_3^6+4F_2^4F_4^2F_5^2+4F_1^4F_6^2F_7^2+
4F_1^4F_4^2F_5^2+4F_0^4F_6^2F_7^2+4F_3^4F_6^2F_7^2+4F_0^4F_4^2F_5^2-
2F_0^4F_4^4-2F_1^4F_4^4+4F_3^4F_4^2F_5^2+24F_0^2F_1^2F_2^2F_3^2+
8F_0^2F_2^2F_7^2F_5^2-8F_0^2F_2^2F_7^2F_4^2+8F_1^2F_3^2F_7^2F_5^2-
8F_1^2F_3^2F_7^2F_4^2-8F_1^2F_3^2F_6^2F_5^2+8F_4^2F_6^2F_2^2F_0^2+
8F_4^2F_6^2F_1^2F_3^2-\break 8F_4^2F_6^2F_1^2F_2^2+8F_1^2F_2^2F_4^2F_7^2-
8F_0^2F_3^2F_5^2F_7^2+8F_0^2F_3^2F_5^2F_6^2+8F_0^2F_3^2F_4^2F_7^2-\break
8F_0^2F_3^2F_4^2F_6^2+8F_1^2F_2^2F_5^2F_6^2-8F_1^2F_2^2F_5^2F_7^2+
24F_0^2F_1^2F_4^2F_5^2+24F_4^2F_5^2F_7^2F_6^2+\break 24F_6^2F_7^2F_1^2F_0^2+
24F_4^2F_5^2F_2^2F_3^2$}
\proclaim
{Proposition}
{The ring
$$\cz[F_a]/(P,Q)$$
is normal of Krull dimension $6$.}
RiN%
\finishproclaim
{\it Proof.\/}
Since it is an complete intersection, one has only to show the the
codimension of the singular locus is at least two. This can be checked
with a computer.\qed
\smallskip
We need the invariant ring
$$B(\Gamma_3[2]):=\bigl(\cz[F_a]/(P,Q)\bigr)^{N_3'}=
\cz[F_a]^{N_3'}/(P,Q).$$
The zero locus of $\vartheta[0]$ is an irreducible subset of
$\calH_3/\Gamma_3[2]$. The restriction of the ring of modular forms
$A(\Gamma_3[2])$ to this hyperelliptic component hence is
$$A(\Gamma_3[2]/\rad((\vartheta[0]^2))$$
\proclaim
{Theorem}
{The restriction of the ring of modular forms
$A(\Gamma_3[2])$ to the hyperelliptic component, set theoretically defined by
$\vartheta[0]=0$, is
$$A(\Gamma_3[2])/(\rad((\vartheta[0]^2))\cong B(\Gamma_3[2])=\cz[F_a]^{N_3'}/(P,Q).$$
It is generated by 14 forms of weight $2$ and 14 forms of weight $3$.
The Hilbert function of this ring is the product of the Hilbert function
of $A(\Gamma_3[2])$ (see \RA) with $(1-z)(1-z^4)/(1-z^8)$,
$$1+14z^2+14z^3+105z^4+175z^5+546z^6+946z^7+2057z^8\cdots$$}
RBh%
\finishproclaim
For the proof we have
to describe the relations between the generators.
For this we need a method which allows to decide whether a system
$f_1,\dots,f_m$ of homogenous polynomials in the variables $F_a$
is linear independent in the quotient
$\cz[F_a]/(P,Q)$ . Of course we assume that all $f_i$ have the same degree.
We used the following method:
First of all, we computed a Groebner basis of the ideal $(P,Q)$, using the
computer algebra {\dun SINGULAR}. Then we computed  the {\it normal forms\/}
$g_1,\dots,g_m$
of the $f_i$ with respect to this basis. The point now is that the
$f_1,\dots,f_m$ are linearly independent in the factor ring $\cz[F_a]/(P,Q)$
if and only if the $g_1,\dots,g_m$ are linearly independent in the polynomial
ring $\cz[F_a]$. This is of course a problem in linear algebra. Using this method
we proved:
\proclaim
{Proposition}
{The space of forms in $B(\Gamma_3[2])$ of weight four,
which is generated by products of two forms of weight
two has dimension\/ $105$. The space of forms of weight six,
which is generated by products of three forms of weight
two has dimension\/ $546$.
\smallni
{\bf Corollary 1.} The\/ $14$ forms of weight two in generate in $B(\Gamma_3[2])$ all
forms of weight\/ $4$ and\/ $6$.
\smallni
{\bf Corollary 2.} There are\/ $560$ monomials of degree\/ $4$ in the\/ $14$ generators
of weight two. Hence there is a\/ $14$-dimensional space of cubic relations.
}
Fwe%
\finishproclaim
It seems to be natural now to consider the subring
$$B'(\Gamma_3[2])\subset B(\Gamma_3[2])$$
generated by the 14 elements of weight 2. Because of \Fwe\ one can conjecture
that both rings agree in even weights. We will see that this is true (\TRe).
\smallskip
Several times a 14-dimensional space occurred. This depends on the fact that
the group $\O(6,\fz_2)\cong S_8$ admits an irreducible irreducible
14-dimensional representation.
There are actually to isomorphism classes which are intertwined by the signum
character.
All 14-dimensional spaces, which occurred so far, are
irreducible.
Next we describe the 14 cubic relations:
\proclaim
{Lemma}
{Restricting to the hyperelliptic component ,the following ternary polynomial is contained in the ideal generated by the
quartic Riemann relations:\smallni
$
18\cdot \Theta_6\cdot \Theta_{11}\cdot \Theta_{14}+18\cdot \Theta_8\cdot \Theta_9\cdot
\Theta_{14}+18\cdot \Theta_2\cdot \Theta_{13}\cdot \Theta_{14}+9\cdot \Theta_1\cdot
\Theta_{12}^2+9\cdot \Theta_1\cdot \Theta_2^2+18\cdot \Theta_3\cdot
\Theta_{14}\cdot \Theta_{15}-43\cdot \Theta_1^3+9\cdot \Theta_1\cdot \Theta_9^2+2\cdot
\Theta_4^3+2\cdot \Theta_8^3+2\cdot \Theta_{11}^3+2\cdot \Theta_{13}^3+18\cdot \Theta_2
\cdot \Theta_5\cdot \Theta_6-18\cdot \Theta_5^2\cdot \Theta_9
+2\cdot \Theta_3^3+9\cdot \Theta_1\cdot \Theta_5^2-18\cdot \Theta_7\cdot \Theta_{11}^2-18
\cdot \Theta_4\cdot \Theta_{14}^2+2\cdot \Theta_{12}^3+18\cdot \Theta_4\cdot \Theta_5\cdot
\Theta_8-18\cdot \Theta_2\cdot \Theta_{15}^2+9\cdot \Theta_1
\cdot \Theta_3^2+18\cdot \Theta_6\cdot \Theta_{12}\cdot \Theta_{13}+18\cdot \Theta_8\cdot
\Theta_{11}\cdot \Theta_{15}+9\cdot \Theta_1\cdot \Theta_{15}^2+18\cdot \Theta_7\cdot
\Theta_{10}\cdot \Theta_{14}+9\cdot \Theta_1\cdot \Theta_8^2+2\cdot \Theta_{14}^3+
18\cdot \Theta_2\cdot \Theta_7\cdot \Theta_8+2\cdot \Theta_9^3+9\cdot \Theta_1\cdot
\Theta_{10}^2+18\cdot \Theta_3\cdot \Theta_{10}\cdot \Theta_{12}+9\cdot \Theta_1\cdot
\Theta_4^2+2\cdot \Theta_5^3+18\cdot \Theta_3\cdot \Theta_6\cdot \Theta_8-18\cdot \Theta_5
\cdot \Theta_9^2+2\cdot \Theta_{15}^3+2\cdot \Theta_6^3+18\cdot \Theta_7\cdot \Theta_9\cdot
\Theta_{13}+18\cdot \Theta_7\cdot \Theta_{12}\cdot \Theta_{15}+30\cdot \Theta_1\cdot
\Theta_7\cdot \Theta_{11}+30\cdot \Theta_1\cdot \Theta_3\cdot \Theta_{13}+18\cdot
\Theta_2\cdot \Theta_{11}\cdot \Theta_{12}+9\cdot \Theta_1\cdot \Theta_{13}^2+30\cdot
\Theta_1\cdot \Theta_4\cdot \Theta_{14}+30\cdot \Theta_1\cdot \Theta_2\cdot \Theta_{15}+18\cdot
\Theta_5\cdot \Theta_{10}\cdot \Theta_{15}+9\cdot \Theta_1\cdot \Theta_{14}^2+18\cdot \Theta_4
\cdot \Theta_9\cdot \Theta_{12}-6\cdot \Theta_{13}\cdot \Theta_{14}\cdot \Theta_{15}-6\cdot
\Theta_5\cdot \Theta_8\cdot \Theta_{14}-6\cdot \Theta_3\cdot \Theta_4\cdot
\Theta_{15}-6\cdot \Theta_9\cdot \Theta_{10}\cdot \Theta_{15}-6\cdot \Theta_6\cdot
\Theta_7\cdot \Theta_{14}-6\cdot \Theta_4\cdot \Theta_5
\cdot \Theta_{12}-6\cdot \Theta_7\cdot \Theta_8\cdot \Theta_{15}-6\cdot \Theta_3\cdot
\Theta_7\cdot \Theta_9-6\cdot \Theta_9\cdot \Theta_{12}\cdot \Theta_{14}-6\cdot \Theta_4\cdot
\Theta_8\cdot \Theta_9-6\cdot \Theta_2\cdot \Theta_7\cdot \Theta_{12}-6\cdot \Theta_3\cdot
\Theta_6\cdot \Theta_{12}-6\cdot
\Theta_4\cdot \Theta_7\cdot \Theta_{10}-6\cdot \Theta_{11}\cdot \Theta_{12}\cdot
\Theta_{15}-6\cdot \Theta_6\cdot \Theta_8\cdot \Theta_{13}-6\cdot \Theta_5\cdot
\Theta_7\cdot \Theta_{13}-6\cdot \Theta_{10}\cdot \Theta_{12}\cdot \Theta_{13}+9\cdot
\Theta_1\cdot \Theta_6^2-6\cdot \Theta_3\cdot
\Theta_8\cdot \Theta_{10}-6\cdot \Theta_9\cdot \Theta_{11}\cdot \Theta_{13}-6\cdot
\Theta_2\cdot \Theta_4\cdot \Theta_{13}-6\cdot \Theta_5\cdot \Theta_6\cdot \Theta_{15}-6\cdot
\Theta_2\cdot \Theta_6\cdot \Theta_9-6\cdot \Theta_2\cdot \Theta_5\cdot \Theta_{10}-6\cdot
\Theta_4\cdot \Theta_6\cdot \Theta_{11}
-6\cdot \Theta_{10}\cdot \Theta_{11}\cdot \Theta_{14}-6\cdot \Theta_2\cdot
\Theta_8\cdot \Theta_{11}-6\cdot \Theta_2\cdot \Theta_3\cdot \Theta_{14}+18\cdot
\Theta_5\cdot \Theta_{11}\cdot \Theta_{13}+18\cdot \Theta_2\cdot \Theta_9\cdot
\Theta_{10}+18\cdot \Theta_4\cdot \Theta_6\cdot \Theta_7+
18\cdot \Theta_6\cdot \Theta_9\cdot \Theta_{15}+18\cdot \Theta_4\cdot \Theta_{13}\cdot
\Theta_{15}+18\cdot \Theta_8\cdot \Theta_{10}\cdot \Theta_{13}+18\cdot \Theta_2\cdot
\Theta_3\cdot \Theta_4+18\cdot \Theta_4\cdot \Theta_{10}\cdot \Theta_{11}+9\cdot
\Theta_1\cdot \Theta_{11}^2+
18\cdot \Theta_3\cdot \Theta_9\cdot \Theta_{11}-18\cdot \Theta_2^2\cdot \Theta_{15}-18\cdot
\Theta_6\cdot \Theta_{10}^2-18\cdot \Theta_6^2\cdot \Theta_{10}-18\cdot \Theta_3^2\cdot
\Theta_{13}-18\cdot \Theta_3\cdot \Theta_{13}^2-18\cdot \Theta_8
\cdot \Theta_{12}^2-18\cdot \Theta_8^2\cdot \Theta_{12}-18\cdot \Theta_4^2\cdot
\Theta_{14}+2\cdot \Theta_{10}^3+2\cdot \Theta_7^3+30\cdot \Theta_1\cdot \Theta_8\cdot
\Theta_{12}+30\cdot \Theta_1\cdot \Theta_5\cdot \Theta_9+2\cdot \Theta_2^3+
30\cdot \Theta_1\cdot \Theta_6\cdot \Theta_{10}-6\cdot \Theta_3\cdot \Theta_5\cdot
\Theta_{11}+18\cdot \Theta_3\cdot \Theta_5\cdot \Theta_7-18\cdot \Theta_7^2\cdot
\Theta_{11}+9\cdot \Theta_1\cdot \Theta_7^2+18\cdot \Theta_5\cdot
\Theta_{12}\cdot \Theta_{14}$
\smallni
Its orbit under the group $\Gamma_{3,\theta'}$ gives a system of\/ 14 relations.}
Relt%
\finishproclaim
The relations of degree $\le 3$ above do not generate the ideal of all relations.
We also need some quartic relations. To compute them one needs the expressions
of the $\vartheta[m]^4$ as linear combinations of the $\Theta_1,\dots,\Theta_{15}$.
Here is one:
$$6\vartheta\left[\matrix{0\cr0\cr0}\matrix{0\cr0\cr1}\right]^4=
\Theta_{13}+\Theta_{14}+\Theta_{15}.$$
(We use the notation $\vartheta[m]=\vartheta[m',m'']$.)
The others are obtained by applying $\O(\fz_2^6)$.
Actually there will be two types of quartic relations between the $\vartheta[m]^4$.
Let $M$ be a set of characteristics
We denote by $\vartheta[M]$ the product of all $\vartheta[m]$ with $m\in M$.
Let  $M$ now be a two dimensional $\fz_2$-vector space of characteristics.
There are three cosets $M$, $a+M+a$, $b+M$ containing only even characteristics.
We  have the Riemann relation
$$\vartheta [M]= \vartheta [M+a] \pm  \vartheta [M+b].$$
If we set $\vartheta[0]$ to zero we obtain the relation
$\vartheta [a+M]^4 = \vartheta [b+M]^4$.
\proclaim
{Proposition}
{There are\/ $210$ relations between the $\Theta_i$
(considered on the hyperelliptic component $\vartheta[0](\tau)=0$) which come
from the relations
$$\vartheta [a+M]^4 = \vartheta [b+M]^4.$$
Here $M$ is a two-dimensional $\fz_2$-vector space of characteristics
and $a+M$ and $b+M$ are the two orbits consisting only of even characteristics.}
Qur%
\finishproclaim
One can expand these relations as quartic polynomials in the $\Theta_i$.
We don't print them here.
\medskip
We also need quartic relations which contain only squares of  thetas :
 They are of the form
$$\vartheta[{m_1}]^2\vartheta[{m_2}]^2 \pm\vartheta[{m_3}]^2\vartheta[{m_4}]^2=
\vartheta[ {m_5}]^2\vartheta[{m_6}]^2 \pm\vartheta[{m_7}]^2\vartheta[{m_8}]^2.$$
Now assuming
$\vartheta[{m_1}]=\vartheta[{0}]=0$
we have
$$ \pm\vartheta[{m_3}]^2\vartheta[{m_4}]^2=\vartheta[{m_5}]^2
\vartheta[{m_6}]^2 \pm \vartheta[{m_7}]^2\vartheta[{m_8}]^2.$$
Squaring we get
$$ \vartheta[{m_3}]^4\vartheta[{m_4}]^4-\vartheta[{m_5}]^4\vartheta[{m_6}]^4 -
\vartheta[ {m_7}]^4\vartheta[{m_8}]^4=
  \pm2\vartheta[{m_5}]^2\vartheta[{m_6}]^2 \vartheta[{m_7}]^2\vartheta[{m_8}]^2:$$
Squaring again we get a quartic relation  among $\vartheta[m]^4$.
This can be expressed in the $\Theta_i$.
The total number of these relations between the $\Theta_i$ is 105.
\proclaim
{Proposition}
{The relations of the type (considered in $B(\Gamma_3[2])$)
$$ (\vartheta[{m_3}]^4\vartheta[{m_4}]^4-\vartheta[{m_5}]^4\vartheta[{m_6}]^4 -
\vartheta[ {m_7}]^4\vartheta[{m_8}]^4)^2=
  4\vartheta[{m_5}]^4\vartheta[{m_6}]^4 \vartheta[{m_7}]^4\vartheta[{m_8}]^4$$
  give\/ $105$ relations between the $\Theta_i$}.
Thrv%
\finishproclaim
We will not give the explicit polynomials in the $\Theta_i$.
The relations described so far are still not all. The  Schottky relation
$$\Bigl(\sum_m\vartheta[m]^8\Bigr)^2-8\sum_m\vartheta[m]^{16}$$
is an extra relation.
The expansion of this relation as polynomial in the $\Theta_i$ is very big.
We don't print it.
\smallskip
Remember that we now have the following system of relations between the
15 functions $\Theta_i$ considered in $B(\Gamma_3[2])$.
\vskip2mm
\item{1)} One linear relation,
\item{2)} 14 cubic relations,
\item{3)} two systems of quartic relations one consisting of 105 the other of 210
relations,
\item{4)} an extra quartic relation.
\smallni
This set of relations is permuted under the action of $\Gamma_{3,\vartheta}'$.
We consider the ideal
$\calR$ generated by all these relations.
\smallskip
Using the computer algebra system {\dun  SINGULAR} one can get a Gr"obner basis
of this ideal, which allows to do several computations, for example it is possible
to get the Hilbert function of the algebra
$$\cz[T_1,\dots,T_{15}]/\calR.$$
Here the $T_i$ are formal variables which stand for the $\Theta_i$.
\smallskip
{\dun SINGULAR} gives the following Hilbert function:
$$
(1+8z^2+36z^4+106z^6+91z^8+14z^{10})/(1-z^2)^6.
$$
(The weight of $\Theta_i$ is two.)
One checks immediately that this series agrees with the even part of the Hilbert series
of the ring $B(\Gamma_3[2])$ as has been described in 6.2. Hence we obtain that
our ideal describes all relations.
\smallskip
We mention that the our system of relations is not minimal.
The system of relations has the advantage to be
invariant under the group $\Gamma_{3,\vartheta}'$.
\smallskip
The final result is:

\proclaim
{Theorem}
{The subring $B^{(2)}(\Gamma_3[2])$ of forms of even weight of the ring
$$B(\Gamma_3[2])=A(\Gamma_3[2])/(\rad((\vartheta[0]^2))$$
is generated by\/ $15$ forms $\Theta_1,\dots,\Theta_{15}$. The ideal of relations
is generated by one linear relation,\/ $14$ ternary relations and\/ $127$
quartic relations: The Hilbert function is given by the formula
$$\eqalign{&
{1+8z^2+36z^4+106z^6+91z^8+14z^{10}\over (1-z^2)^6}=\cr&
1+14z^2+105z^4+546z^6+2057z^8+6062z^{10}+14945z^{12}+\cr&
32306z^{14}+63217z^{16}+114478z^{18}+\cdots
\cr}$$
The projective variety of this ring is the hyperelliptic component in the Satake
compactified Siegel modular variety $\overline{\calH_3/\Gamma_3[2]}$. Hence this closure
is a normal variety.
}
TRe%
\finishproclaim
\neupara{Blowing up}%
We want to investigate the homomorphism
$$B^{(2)}(\Gamma_3[2])=\cz[\dots\vartheta[m]^4\dots]\lo\cz[Y_1,\dots,Y_{14}]/\calJ.$$
Here the $\vartheta^4[m]$ are understood as forms on the hyperelliptic component.
Recall that  the $\vartheta^4[m]$ can be expressed by means of the $\Theta_i$.
Hence the ideal of relations between the $\vartheta^4[m]$ (restricted to the hyperelliptic
component) is known from \TRe. The right hand side is the graded algebra
of the ball quotient. The ideal $\calJ$ of relations has also been described
explicitly in \LK. Recall that we have an explicit isomorphism $\O(\fz_2^6)\cong S_8$.
The group $\O(\fz_2^6)$ acts transitively on the even characteristics $m\ne0 $.
The group $S_8$ acts linearly on the $Y_i$. (They are polynomials in 8 variables
$X_1,\dots, X_8$ on which the group $S_8$ acts by permutation.)
\smallskip
Using Thomae's theorem \ST\ and the expressions of the $Y_i$ as polynomials in the
$X_i$ (section 1) one can compute explicit expressions of the images of the
$\vartheta[m]^4$. The base locus ideal is the ideal generated by the images.
Since the above homomorphism is $\O(\fz_2^6)=S_8$-equivariant, it is sufficient to give
the image of one
$\vartheta[m]^4$.
\proclaim
{Lemma}
{The homomorphism
$$\cz[\dots\vartheta[m]^4\dots]\lo\cz[Y_1,\dots,Y_{14}]/\calJ$$
is defined by{\ninepoint
$$\vartheta\left[\matrix{0&0\cr0&0\cr0&1\cr}\right]^4\loma
(Y_1-Y_{10}+Y_{11}-Y_{14})(Y_1-Y_2-Y_6+Y_7+Y_8-Y_9+Y_{11}-Y_{13})(Y_8-Y_9).$$}%
The others are obtained applying $S_8$.}
Ioah%
\finishproclaim
The base locus is defined by the ideal which is generated by the
images of the $\vartheta[m]^4$.
\proclaim
{Proposition}
{The ideal, which is generated by the images of the $\vartheta[m]^4$ in the ring
$\cz[Y_1,\dots,Y_{14}]/\calJ$
is the intersection of 56 ideals  which are
generated by linear forms. One of them is
$$(Y_2,Y_3,Y_5,Y_6,Y_8,Y_9,Y_{10},Y_{11},Y_{13},Y_{14})$$
The others are obtained applying $S_8$.}
LFb%
\finishproclaim
This can be verified by a {\dun SINGULAR} calculation. However:
in the last section we shall give a more geometric  and intrinsic description of the components of  zero set of the ideal  generated by the images of the $\vartheta[m]^4$. \bigskip
\noindent
Let $\gotA$ be an ideal in a (commutative and with unity) noetherian ring $A$ .
The blow up of $A$ along $\gotA$ is the graded $A$-algebra
$$\Bl(A,\gotA):=\bigoplus_{n=0}^\infty\gotA^n.$$
We can consider $\Bl(A,\gotA)$ as a subring of $A[T]$, using the embedding
$$\Bl(A,\gotA)\lo A[T],\quad a\in\gotA^n\loma aT^n.$$
Now we assume
$$A=\cz[X_1,\dots,X_n]/\gota.$$
We choose polynomials $Q_1,\dots,Q_m$ in $\cz[X_1,\dots,X_n]$, whose
images in $A$ generate $\gotA$.
We consider the homomorphism of polynomial rings
$$\cz[X_1,\dots,X_n,Y_1,\dots,Y_m]\lo\cz[X_1,\dots,X_n,T],\qquad
Y_i\loma Q_iT.$$
We denote by $\tilde{\gotA}$ the
the inverse image of the ideal generated by
$\gota$.
Then we have
$$\Bl(A,\gotA)=\cz[X_1,\dots,X_n,Y_1,\dots,Y_m]/\tilde{\gotA}.$$
We need a modification of this construction. Assume that $A$ already is a
graded algebra and $\gotA$ a graded ideal.
Then $\Bl(A,\gotA)$ is a bigraded algebra and one can define in an obvious way
the projective variety $\Biproj(\Bl(A,\gotA))$ together with a morphism
$$\Biproj(\Bl(A,\gotA))\lo\proj A.$$
We describe this map in terms of coordinates:
For simplicity we assume that the degree of the $X_i$
are one and that the $Q_i$ all are of the same
degree. The variety
$\Biproj(\Bl(A,\gotA))$ then consists
of all pairs $([x],[y])$,
$[x]\in P^{n-1}(\cz)$, $[y] \in P^{m-1}(\cz)$ such that
$P(x,y)=0$ for all polynomials $P(X,Y)\in\tilde{\gotA}$ wich are homogenous in
$X=(X_1,\dots,X_n)$ and $Y=(Y_1,\dots,Y_m)$.
\smallskip
We want to apply this in the following situation:
We consider the rings
$$B=\cz[Y_1,\dots,Y_{14}]/\calJ\quad \text {and}\quad  A=\cz[\dots\vartheta[m]^4\dots].$$
Here the $\vartheta[m]^4$ are understood to be restricted to the hyperelliptic
component, hence $\vartheta[0]^4=0$. We recall that there is a homomorphism
$$A=\cz[\dots\vartheta[m]^4\dots]\lo B.$$
The ideal which we want to blow up
is generated by the images of the $\vartheta[m]^4$ in $B$. We have shown that
this ideal is generated by 56 linear forms. We denote the blow up
of this ideal by $B^*$ and the associated variety by $Y^*$. This is a model
which lies over the Siegel model $X=\proj(A)$ and over the ball-model $Y=\proj(B)$.

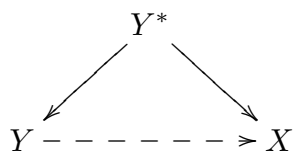
$$\xymatrix{&
  Y^*\ar[dl]\ar[dr]&\\Y\ar@{-->}[rr]&
 &X}$$
We want to get explicit information about the model $Y^*$. Actually it is
too difficult to blow up the 56 components of the base locus in one step.
But this is not necessary. To get information about the blow up of a
small neighborhood of a given point in $Y$, it is sufficient to blow up
those of the linear components which contain this point. To make use of this
we need the intersection behavior of the 56. A direct computation shows:
\proclaim
{Lemma}
{When a subset of the\/ $56$ linear components has a a common intersection,
it intersects also in a cusp. In each cusp meet\/ $8$ of the
linear components.
\smallni
For example the 8 components defined by the following ideals
(each generated by 10 linear forms)
$$\eqalign{
&(Y_2,Y_3,Y_5,Y_6,Y_8,Y_9,Y_{10},Y_{11},Y_{13},Y_{14}),\cr
&(Y_1,Y_2,Y_3-Y_{11}+Y_{14},Y_4-Y_{12}+Y_{14},Y_5,Y_6-Y_{11},Y_7-Y_{12},\cr&
\qquad Y_8-Y_9,Y_{10},-Y_{14}+Y_{13}),\cr
&(Y_1,Y_2,Y_3,Y_5-Y_{13}+Y_{14},Y_6-Y_{13},Y_7-Y_{12},Y_8-Y_{13},\cr&
\qquad Y_9-Y_{14},Y_{10}-Y_{13}+Y_{14},Y_{11}-Y_{13}),\cr
&(Y_1,Y_2-Y_{11}+Y_{13},Y_3-Y_{11}+Y_{13},Y_4-Y_{12},Y_5-Y_{11},Y_6-Y_{11},\cr&
\qquad Y_7-Y_{12},Y_9,Y_{10}-Y_{11},Y_{14}),\cr
&(Y_1,Y_2,Y_3,Y_4-Y_7,Y_5,Y_6,Y_{10},Y_{11},Y_{13},Y_{14}),\cr
&(Y_1-Y_{10}-Y_{14},Y_2-Y_{10}-Y_{14},Y_3,Y_4-Y_9-Y_{12}+Y_{14},Y_5-Y_{10},\cr&
Y_6,Y_7-Y_9-Y_{12}+Y_{14},Y_8,Y_{11},Y_{13}),\cr
&(Y_2-Y_{10},Y_3-Y_{11},Y_4-Y_{12},Y_5-Y_{10},Y_6-Y_{11},Y_7-Y_{12},\cr&
Y_8,Y_9,Y_{13},Y_{14}),\cr
&(Y_1+Y_7-Y_{12},Y_2+Y_{13}-Y_{14},Y_3+Y_{13},Y_4-Y_{12}+Y_{14},\cr&
Y_5,Y_6,Y_8,Y_9,Y_{10},Y_{11}),\cr
}$$
meet in a cusp with the coordinates
$[0,0,0,1,0,0,1,0,0,0,0,1,0,0]$.
}
LiLm%
\finishproclaim
The intersection of the 8 ideals, which are generated by these linear forms
is generated by 17 elements of degree $\le 2$. We will not print
them and mention just that this ideal is simple enough to be blown up
by means of {\dun SINGULAR}.
In this way we obtain for each cusp $s$ a partial blow up $Y_s^*$.
They exhaust $Y^*$ in the following sense. Consider for each
cups  the complement
of the remaining $48=56-8$ components in $Y^*_s$.
These are quasi projective varieties which give an open covering
of $Y^*$. This means that we have a description of $Y^*$ be explicit
equations. We will not print them here but mention just some
consequences which can taken from them.
\proclaim
{Proposition}
{The model $Y^*$ is smooth outside the inverse images of the cusps.
It is not smooth everywhere. The inverse image of a cusp is a two
dimensional irreducible variety.}
PYns%
\finishproclaim
To desingularize one needs one further blow up:
\proclaim
{Proposition}
{The blow up of $Y^*$ along the inverse images of the cusps
is a non-singular model $\tilde Y$.}
nsY%
\finishproclaim
For the proof it is convenient to proceed slightly different.
We first blow up the cusps and then the inverse images of the
56 linear spaces. The result is the same. Direct computations shows:
\proclaim
{Proposition}
{The blow up of the ball model $Y$ along the cusps is a
smooth model $Y'$.
The inverse images of the cusps are irreducible.
The intersection of three different of the\/
strict transforms of
the $56$ linear subspaces of $Y$ is empty.}
Sac%
\finishproclaim
Since $Y'$ is smooth, the inverse image of the union of the
56 linear spaces and the union of the strict transformed only differ
by an invertible ideal. Since the blow-up of two ideals which
differ only by an invertible ideal is the same, we now only have to
blow-up the union of the 56 strict transforms.
For the proof of \nsY\ it now is sufficient to show that the blow up
of $Y'$ along the union of two of the 56 strict transformed is smooth.
Actually the strict transformed are smooth and the local analytical
behavior in an intersection point is the same as $z_1=z_2=$ and
$z_3=z_4=$ in a $\cz^5$. This shows that the blow-up remains smooth
and proves \nsY.
\smallskip
Finally we compare the model $\tilde Y$ with $\bar M_{0,8}$.
Recall that $M_{0,8}$ is the moduli space of 8 ordered points on a
projective line and $\bar M_{0,8}$ the Mumford-compactification
by marked stable curves.
It is known [Ka] (s.~also [AL]) that this is a projective smooth variety
and that there exist regular contraction maps ([Ka], s.~also [AL])
$$\bar M_{0,8}\lo Y.$$
From the universal property of blowing up one obtains a regular map
$$\bar M_{0,8}\lo \tilde Y.$$
We claim that this is biholomorphic.
\proclaim
{Theorem}
{The two diagrams
$$\xymatrix{&
  \tilde Y\ar[dl]\ar[dr]&\\Y\ar@{-->}[rr]&
 &X}\qquad
\xymatrix{&
  \bar M_{0,8}\ar[dl]\ar[dr]&\\Y\ar@{-->}[rr]&
 &X}
 $$
are isomorphic.
}
HS%
\finishproclaim
The proof rests on the following simple

\proclaim
{Lemma}
{Let $V,W$ be smooth complete varieties and $f:V\to W$
a regular map which induces a biholomorphic map $V_0\to W_0$
of certain Zariski
open subsets $V_0\subset V$, $W_0\subset W$. Assume that $V-V_0$ and
$W-W_0$ are of pure codimension one and that the number of irreducible
components in both cases is the same. Then $f$ is biholomorphic.}
Lxy%
\finishproclaim
We omit the simple proof. We have to apply it to
$V=\bar M_{0,8}$, $V_0=M_{0,8}$, $W=\tilde Y$ and $W_0$ the complement of the
inverse images of the
Heegner divisors in $Y$.
\smallni
{\it Proof of \HS\/.}
Firstly we recall shortly the irreducible components of $\bar M_{0,8}-M_{0,8}$.
Their generic points are marked curves $(C,x_1,\dots,x_8)$, where
$C$ is the union of two
$P^1$ with one intersection point (normal crossing). We have to choose $a$
point on the first $P^1$ and $b$ on the second one, $a+b=8$. The possibilities
to get a stable marked curve are $(a,b)=(6,2),\;(5,3),\; (4,4)$.
In the case $(a,b)=(6,2)$ we have $28={8\choose 2}$ possibilities, in the case
$(a,b)=(5,3)$ there are $56={8\choose 3}$ and in the case
$(a,b)=(4,4)$ there are $35={8\choose 4}/2$ possibilities. The denominator
$2$ comes from the fact that the role of the two $P^1$ can be interchanged.
Hence we have $28,56$ and $35$  irreducible components. The group
$S_8$ acts on each of the three systems transitively.
\smallskip
This picture can be recovered in $\tilde Y$. We have 28 irreducible subvarieties
of codimension 1 which map to the 28 Heegner divisors in $Y$. Hence these
28 are visible already in $Y$.
Then we have the 56 irreducible subvarieties
of codimension 1 which come from blowing up the 56 components of the
base locus. These 56 are visible in $Y^*$.
Finally we have the 35 irreducible subvarieties
of codimension 1 which come from blowing up the cusps.
\smallskip
Since $\O(\fz_2^6)$ permutes the three systems, this counting also shows:
\proclaim
{Proposition}
{The boundary $\bar M_{0,8}-M_{0,8}$ is the union of\
$28+56+35$ irreducible subvarieties of codimension one which correspond to
distributions $(6,2)$, $(5,3)$, $(4,4)$ of eight points on two
crossing $P^1$. Under the map $\bar M_{0,8}\to Y$
type $(6,2)$-components map to  the $28$
Heegner divisors, type $(5,3)$-components to the
$56$ components of the base locus and type $(4,4)$-components are contracted
to the cusps.}
Ptet%
\finishproclaim
\smallskip
In the next section we shall consider also the image  of  the   boundary $\bar M_{0,8}-M_{0,8}$
under the map $\bar M_{0,8}\to X$
\neupara{A combinatorial  approach}%
In this section we want to give a combinatorial  approach to the
two different compactifications  $X$ and $Y$ of the hyperelliptic component.
We will restrict our attention to those points in $X$ and $Y$ that do
not  correspond to non singular complete hyperelliptic curves.
We  shall denote by $X_0$ and $Y_0$ the open set in $X$ and $Y$
corresponding to hyperelliptic curves.
Thus we need at least a set theoretical description of the complementary  loci.
\smallskip
Firstly we recall shortly the structure of $X- X_0$.
We use the embedding of the Siegel half plane $\calH_g$ into a Grassmanian and
we denote by $\calH_g^*$ the union of $\calH_g$ with the rational boundary
components. The Satake compactification then is $\calH_g^*/\Gamma_g[q]$.
Details can be found for example in [Fr3].
\smallskip
We restrict now to $g=3$.
Recall that $X_0$ is the subset of all points $\calH_3/\Gamma_3[2]$, which are
represented
by a (smooth) hyperelliptic Riemann surface for which $\vartheta[0]$ vanishes
and that $X$ is the closure of $X_0$ in the Satake compactification.
We recall that a point $Z\in\calH_3$ is reducible if
it is conjugate  with respect to the action of $\Gamma_3$ to a point of the  form
$$\pmatrix{\tau_1& 0\cr 0&\tau_2}$$
with $\tau_i\in \calH_{g_i}, \,\,\, g_1+g_2=3$.
We shall denote by  $\calR_3$ the set of reducible points, on which $\vartheta[0]$
vanishes. It is known that  $\calR_3/\Gamma_3[2]$
is contained in $X$ and even more
$$X=X_0\cup \overline{\calR_3/\Gamma_3[2]}\qquad\hbox{(disjoint union)}.$$
We want to count how many  reducible components   appear in $X-X_0$.
Each reducible component
of $\overline{\calR_3/\Gamma_3[2]}$
is defined
by the vanishing of 6 thetanullwerte, (one of them zero, since we
work on the locus $\vartheta[0]=0$). We refer to [Gl] for details.
A sextuplet of characteristics $m_1, \dots,m_6$ corresponding to these six
thetanullwerte have the property that all sums of three of them result
to be an odd characteristic. An example is given by
\bigni
\centerline{\vbox{
\halign{$\hfil\ #\ \hfil$&$\hfil\ #\ \hfil$&$\hfil\ #\ \hfil$
&$\hfil\ #\ \hfil$&$\hfil\ #\ \hfil$
&$\hfil\ #\ \hfil$\cr
m_1&m_2&m_3&m_4&m_5&m_6\cr
\noalign{\vskip3pt}
0&1&0&1&0&0\cr 0&1&0&0&1&0\cr 0&1&0&0&0&1\cr
\noalign{\vskip1pt}
0&1&1&0&1&1\cr 0&0&1&0&0&1\cr0&1&1&0&0&0\cr}}}
\medni
Hence we obtain:
\proclaim
{Proposition}
{The boundary part $X-X_0$ of the hyperelliptic component $X$ consists of\/
$56$ irreducible 4-dimensional components. They correspond
to sets of 6 even characteristics, one of them $0$, such that all sums of three
of them are odd. The component is the intersection of the corrseponding
zero divisors $\vartheta[m]=0$. The group $\Gamma_{3,\vartheta}$ permutes
the components transitively.
}
Pfitc%
\finishproclaim
Part of the boundary $X-X_0$ is the intersection of $X$ with the Satake boundary.
We want to describe this part more closely. Recall that
$$\calH_3^*-\calH_3$$
is the disjoint union of the $\Gamma_3$-orbits
of the three standard components
$$\pmatrix{\imag\infty&0&0\cr0&\imag\infty&0\cr 0&0&\imag\infty\cr},\qquad
\pmatrix{\tau&0&0\cr0&\imag\infty&0\cr 0&0&\imag\infty\cr}\ (\tau\in\calH_1),\qquad
\pmatrix{\tau& 0\cr 0&\imag\infty}\ (\tau\in\calH_2).$$
We denote by $\calD_3$ the part of $\calH_3^*-\calH_3$ on which $\vartheta[0]$
vanishes.
We see that  $\calD_3/\Gamma_3[2]$ is contained already
in the closure of the reducible locus. Hence we have the inclusions
$$\calD_3/\Gamma_3[2]\subset \overline{\calR_3/\Gamma_3[2]}
\subset X\subset \overline{\calH_3/\Gamma_3[2]}.$$
\smallskip
Next we want to describe the zero dimensional boundary components of $X$.
These are the points which come from the orbit of the point
$$\pmatrix{\imag\infty&0&0\cr0&\imag\infty&0\cr0&0&\imag\infty\cr}.$$
 In [Gl] one finds the description of the zero dimensional boundary components
as intersections of 28 zero divisors $\vartheta[m]=0$.
This description can be related to stars. Recall that a star is a set of
$4$ odd characteristics with a certain property (\DStar).
\proclaim
{Proposition}
{Let $S$ be a star.
 There are 27  even  characteristics  $m\ne 0$
 orthogonal  to at least one  characteristic in $S$.
 The intersection of the 27 divisors $\vartheta[m]=0$ is a zero-dimensional
 boundary component of $X$. This gives a 1-1-correspondence between the
 105 stars and the 105 zero-dimensional boundary components.}
 PSiX%
 \finishproclaim
 We mention also how the remaining 8 even characteristics can be described in a nice
 way:
 The star $S$ spans a three dimensional space  M.
 There exists a unique even characteristic $n$ such that
 $M+n$ contains only even
 characteristics. The entries of this coset are the 8
 remaining characteristics.
\smallskip
Now we describe the closures of the three dimensional boundary components of $X$. They
are of the type
$$\pmatrix{\tau&0\cr 0&\imag\infty\cr},
\quad\tau\in\calH_{2}.$$
In [Gl] one finds also  the description of the three dimensional boundary components
as intersections of 16 zero divisors $\vartheta[m]=0$. The description can be related to the odd characteristics

\proclaim
{Proposition}
{Let $\alpha$ be an odd characteristic,
there are exactly 15 even  characteristics  $m\ne 0$ orthogonal to it.
The corresponding divisors $\vartheta[m]=0$ intersect in a
three-dimensional boundary component of $X$. This gives a  a 1-1-correspondence between the
 28 odd characteristics  and the  $28$ irreducible components of the boundary $\calD_3 /\Gamma_3[2]$ which are all three dimensional.
}
odbC%
\finishproclaim
A similar computation can be done for the one dimensional
boundary components of $\overline{\calH_3/\Gamma_3[2]}$,
knowing that in this case each component is given by the vanishing of
24  suitable thetanullwerte, cf [Gl].
In this case we get:

 \proclaim
{Proposition}
{Let $\alpha_1$ and $\alpha_2$  be   odd characteristics such that $\alpha_1+\alpha_2$
is an even characteristics, then there are   exactly 23  even
characteristics $m\ne 0$ orthogonal to $\alpha_1$ or to $\alpha_2$ .
The corresponding divisors $\vartheta[m]=0$ intersect in a
one-dimensional boundary component of $X$. This gives a  a 1-1-correspondence
between the
 $210$ pairs of  odd characteristics  of the above form  and the  $210$ one
 components of the boundary $\calD_3 /\Gamma_3[2]$. }
Pbca%
\finishproclaim
Next we recall the structure of $Y-Y_0$. We have two possibilities, to use
the {\dun GIT}-picture or the ball-picture. The {\dun GIT}-picture has been
explained in
[Koi]. We just recall  the following: Consider the map
$$\cz^8\lo (P^1)^8\lo (P^1)^8//\SL(2)=Y.$$
In section one we used the variables $X_1,\dots,X_8$ to describe this $\cz^8$.
We have to consider the 28 differences $W_{ij}=X_i-X_j$, $1\le i<j\le 8$.
The images of their zero set in $Y$ describe $28$ irreducible divisors.
There union is the boundary part $Y-Y_0$. We now switch to the ball picture.
The reason is that in this picture the quadratic space $\fz_2^6$ occurs, whose
elements also play the role of theta characteristics.
\smallskip
In \RSf\ we showed that
the 28 pairs $(i,j)$ are in 1-1-correspondence to the anisotropic elements of
$\fz_2^6$, in other words to the odd characteristics.
Recall that to each such element $\alpha$ we introduced an irreducible
divisor (Heegner divisor)  $\calB_\alpha$ in $\overline{\calB/\Gamma[1+\imag]}$.
\proclaim
{Proposition}
{In our identification $Y=\overline{\calB/\Gamma[1+\imag]}$
the boundary divisor in $Y$ defined by ``$W_{ij}=0$", cf.~Proposition \VExH,
corresponds to the
Heegner divisor $\calB_\alpha$, where $\alpha$ is the anistropic
element of $\fz_2^6$ related to $(i,j)$ as in \RSf.}
PiWij%
\finishproclaim
For a sequence of  odd characteristics
$\alpha_1, \dots, \alpha_k $,  we associate
the intersection
$$\calB_{\alpha_1,\dots, \alpha_k}:=\calB_{\alpha_1}\cap\dots
\cap\calB_{\alpha_k} .$$
Now we want to compare the Siegel- and the ball picture, i.e.~we have to consider
the rational map $\xymatrix{Y\ar@{-->}[r] &X\\}$.
 \proclaim
{\bf Proposition}
{Under the rational map
$\xymatrix{Y\ar@{-->}[r]&X\\}$ the 28 Heegner divisors are contracted to the
28 three-dimensional boundary components. The  210 intersections
$\calB_{\alpha_1, \alpha_2}$   with $\alpha_1 +\alpha_2$ an even characteristic
are contracted to the 210 one dimensional boundary components.
The 420 intersections
$\calB_{\alpha_1,\alpha_2,\alpha_3} $ with $\alpha_i +\alpha_j$
even characteristics ( $1\leq i<j\leq 3$) are contracted to the 105
zero dimensional boundary components. More precisely,  the  four
$\calB_{\alpha_1,\alpha_2,\alpha_3} $ obtained by a star map to the same point .}
Urmc%
\finishproclaim

 The rational map  is defined outside the base locus,
which has been
described in \LFb\ as a concrete three dimensional variety. When $B$ is a divisor in $X$
we can take its inverse image in the complement of the base locus and then take the
closure $A$ in $X$. We call $A$ the pull-back of $B$.
From \VExH\ and \RSf\ we get:
\proclaim
{Proposition}
{The pull-back of the zero divisor of the $\vartheta[m]$, $m\ne 0$,
is
the union of\/ $12$ Heegner divisors $\calB_\alpha$,
namely those such that $\alpha$ is orthogonal to $m$.}
ZdcX%
\finishproclaim
Recall (\DStar) that stars consist of 4 anisotropic vectors with a certain property.
The union of the four $\calB_\alpha$ is the star-divisor $\calB_S$. A direct inspection
shows (one $m$ is enough):
\proclaim
{Remark}
{Let $m\ne 0$ be an even characteristic. There are five different ways to write the
set of all anistropic $\alpha$, which are orthogonal to $m$, as union of
three stars. The pull-back of the zero divisor of $\vartheta[m]$ is the
union of  any of such three star divisors.}
CcSd%
\finishproclaim
We come now to a description of the cusps. It is known that there are
35 cusps (elements of $\overline{\calB/\Gamma[1-\imag]}-\calB/\Gamma[1-\imag]$).
In the {\dun GIT}-model they correspond to the instable points.
From [Koi] or as an almost direct consequence of the computational Lemma \LiLm,
we obtain:
\proclaim
{Remark}
{The intersection of the twelve Heegner divisors, which correspond to
an even non-zero characteristic (\ZdcX), is a cusp.
This gives a bijection between the cusps and the non-zero even characteristics.}
tHDC%
\finishproclaim
Now we will give a combinatorial
description of the base locus of the rational map
$\xymatrix{Y\ar@{-->}[r]& X\\}$.  We recall that
it is  is the union  of 56 linear spaces.
This has been proved in \LFb\ by a computer calculation
(in the strong ideal theoretic sense).
Here will give an instrinsic (set-theoretical)
description of this locus in terms of Heegner divisors.
\smallskip
First of all we need   some preliminary facts.
\proclaim
{Lemma}
{Let
$\alpha_1$ and $\alpha_2$ be two odd characteristics
such that
$$\alpha_3=\alpha_1+\alpha_2$$
is  still an
odd characteristic. Then
each even characteristic  $m\neq 0$ is orthogonal to  at
least one of the three  $\alpha_1, \,\alpha_2,\, \alpha_3. $}
Rab%
\finishproclaim
In fact each odd characteristics $\alpha_1$ is orthogonal to 15
even characteristic $m\ne 0$. Moreover
 25 even characteristics $m\ne 0$  are  orthogonal to at least one of
 two odd characteristics $\alpha_1$ and $\alpha_2$ whose sum is
 still odd, thus both  are orthogonal to 5 even characteristics $m\ne 0$. Obviously
$\alpha_3=\alpha_1 + \alpha_2$ is  still orthogonal to  the
same 5 even characteristics $m$, and there are 10 more that
are orthogonal to $\alpha_3$, but neither orthogonal   to
$\alpha_1$ nor to  $\alpha_2$.
So we get
that  each  even characteristic $m$  is  orthogonal
to at least  one among $ \alpha_1 , \alpha_2$  and $\alpha_3$.
\smallskip
This discussion also shows:
 \proclaim
{Corollary}
{There are exactly 5 even  characteristics
$m\neq 0$   orthogonal to   $\alpha_1, \,\alpha_2,\, \alpha_3$ in \Rab.}
Sab%
\finishproclaim
Let us  consider a different configuration. We assume
that $\alpha_1,\alpha_2$ are odd but
$\alpha_1 + \alpha_2$  is an even characteristic. In this case
23 even non zero  characteristics   are  orthogonal to at least one
of $\alpha_1,\alpha_2$, and both together are orthogonal to 7 even characteristics.
Let us take an odd $ \alpha_3$ such that  $\alpha_1+\alpha_3$ and
$\alpha_2+\alpha_3$ are even characteristics. Then
27 non zero even  characteristics  are  orthogonal to at least one
among   $\alpha_1, \,\alpha_2, \, \alpha_3$.
There is one and only one odd characteristic $\alpha_4$ which added to
$\alpha_1,\alpha_2,\alpha_3$ gives even characteristics, namely
$\alpha_4=\alpha_1+\alpha_2+ \alpha_3$.
Hence the four odd characteristics $\alpha_1,\alpha_2,\alpha_3,\alpha_4$
form  a star. From this discussion follows:
\proclaim
{Lemma}
{Let $\alpha_1, \dots,  \alpha_k $ be a sequence of odd characteristics
such that each even characteristic is orthogonal at  least to a $\alpha_i$,
Then the sequence contains three characteristics, say
$\alpha_1, \,\alpha_2,\, \alpha_3 $, with  $\alpha_3=\alpha_1 + \alpha_2$.}
Tab%
\finishproclaim
For three odd characteristics
$\alpha_1, \,\alpha_2,\, \alpha_3 $
with the property $\alpha_1 + \alpha_2+\alpha_3=0$  we want to consider the intersection
$\calB_{\alpha_1, \alpha_2, \alpha_3}$. One can check that
it coincides already with the intersection of any
two  of the three Heegner  divisors, i.e
$$\calB_{\alpha_1, \alpha_2, \alpha_3}=
\calB_{\alpha_1, \alpha_2}=\calB_{\alpha_1, \alpha_3}=\calB_{\alpha_2, \alpha_3.}$$
A direct computation tells us that there are exactly 56
triplets of the form  $\alpha_1, \,\alpha_2,\, \alpha_3 $,
with  $\alpha_1 +\alpha_2+\alpha_3=0$. With these notations we have the following
\proclaim
{Proposition }
{The base locus of the  rational map from $Y$ to $X$ is  the set
$$V=\bigcup_{\alpha_1+\alpha_2+\alpha_3=0} \calB_{\alpha_1, \alpha_2, \alpha_3}.$$
To each  $\calB_{\alpha_1, \alpha_2, \alpha_3}$  corresponds to a linear
space of Proposition \LFb. }
Uab%
\finishproclaim

The proof is  an immediate consequence of the previous lemmata
and of the fact that the base locus is necessarily  the union of
intersections of Heegner divisors.\smallskip
 We know from the results of the  previous sections that the blow
 up of the  base locus produces 56 divisors in $Y^*$ or in
 $\bar M_{0,8}$. We are interested in the image of these divisors in $X$.
 From the previous  corollary \Sab, we have exactly 5 even
  characteristics orthogonal to $\alpha_1, \,\alpha_2,\, \alpha_3 $.
  Thus the image of $\calB_{\alpha_1, \alpha_2, \alpha_3}$ is defined
  by the vanishing of 6 thetanullwerte whose characteristics are
  $\Gamma_{3,\vartheta}$-conjugate to those at the begin of the section.
\smallskip
Thus, according to our previous result,  we have that they
 define an irreducible component of
 $ \overline{\calR_3/\Gamma_3[2]}$. Hence we get:
\proclaim{Proposition }
{ The blow up of the base locus maps on the reducible locus of the
hyperelliptic modular variety. In particular the blow
up of $\calB_{\alpha_1, \alpha_2, \alpha_3}$ maps onto the
component defined by the vanishing of  the
thetanullwerte $\vartheta[m]$  with $m$ orthogonal to $\alpha_1, \alpha_2, \alpha_3$. }
Vab%
\finishproclaim
Let us conclude considering  the image of the blow up of the cusps.
First of all we observe that
for each even, non zero characteristic $m$ there are $16$ even
characteristics $n_1,\dots, n_{16}$
such that  for each index $i=1, \dots, 16$, $m+n_i=\beta_i$ is an
odd characteristic. These sixteen odd characteristics appear in the
situation that we are going to describe.
\smallskip
We recall from Lemma \LiLm\  and Proposition \Uab\ that
 8 linear components $\calB_{\alpha_1, \alpha_2, \alpha_3}$  intersect in a cusp $m$.
Each triplet is characterized by being orthogonal to
 the characteristic $ m$.
Viceversa, according to Corollary  \Sab\ ,
to each triplet  $\alpha_1, \alpha_2, \alpha_3$ we can associate,
besides the  characteristics
$m$,   four  more even characteristics  that are orthogonal to them.
The union of all these  even characteristics appearing at least once  for the eight
possible triplets
$\alpha_1, \alpha_2, \alpha_3$, related to the even characteristic $m$,  gives
the sixteen  even  characteristics $$n_1,\dots, n_{16}.$$
In [Gl] one finds also  the description of the two dimensional   components conjugate
to $$\pmatrix{\tau_1&0&0\cr0&\tau_2&0\cr 0&0&\imag\infty\cr}\ (\tau_i\in\calH_1)$$
as intersections of 18 zero divisors $\vartheta[n]=0$. The description can be
related to the  characteristics $0,m,n_1\dots, n_{16}$. In fact we  have
\proclaim{Proposition }
{ The blow up of the cusps  maps on   two dimensional   components in the
hyperelliptic modular variety. In particular the blow up of the cusp $m$
maps onto the
component defined by the vanishing of  the
thetanullwerte $\vartheta[n]$, $n=0,m,n_1\dots, n_{16}$
with  $m+n_i$ odd.}
Zab%
\finishproclaim
We remind that
in Proposition \PYns\  is named the two dimensional variety that is the
inverse image of the cusps in $Y^*$.
\bigskip\noindent
{\paragratit References}\medni
\medskip
\item{[AF]} Allcock,\ D.\ and Freitag,\ E.:
{\it Cubic surfaces and Borcherds products},
Commentarii Math.\ Helv. Vol.~{\bf 77}, Issue 2, 270--296 (2002)
\medskip
\item{[AL]} Avritzer, D. Lange, H.:
{\it The moduli spaces of hyperelliptic curves and binary forms\/},
Math. Z. {\bf 242}, 615-632 (2002)
\medskip
\item{[Bo]} Borcherds,\ R.: {\it Automorphic forms
with singularities on Grassmannians,}
Invent. math. {\bf 132}, 491--562 (1998)
\medskip
\item{[Fr1]} Freitag, E.: {\it Some modular
forms related to cubic surfaces\/}
Kyungpook Math. J. {\bf 43}, No.3, 433-462 (2003)
\medskip
\item{[Fr2]} Freitag, E.:
{\it Comparison of different models of the moduli space
of marked cubic surfaces,\/}
Proceedings of Japanese-German Seminar, Ryushi-do, edited by T.~Ibukyama
and W.~Kohnen, 74-79 (2002)
\medskip
\item{[Fr3]} Freitag, E.: {\it Siegelsche Modulfunktionen\/}
Grundlehren der mathematischen Wissenschaften, {\bf 254}
Berlin-Heidelberg-New York: Springer-Verlag (1983)
\medskip
\item{[FS]} Freitag, E. Salvati-Manni, R.:
{\it Modular forms for the
even unimodular
lattice of signature  (2,10),\/} J.~Algebraic Geom. {\bf16}, 753--791 (2007)
\medskip
\item{[Gl]} Glass, J.  :
{\it Theta constants of genus three,\/}
Compos. Math. {\bf 40}, 123-137 (1980).
\medskip
\item{[Ho]} Howe, R.: {\it The classical groups and invariants of
bilinear forms\/}, The Mathematical Heritage of Hermann Weyl
(Durham, NC, 1987), 133--166, Proc. Sympos. Pure Math.
{\bf 48}, Amer. Math. Soc., Providence, RI (1988)
\medskip
\item{[HM1]} Howard, B.J. Millson, J. Snowden, A. Vakil, R.:
{\it The projective invariants of ordered points on the line},
ArXiv Mathematics
e-prints math.AG/0505096  (2007)
\medskip
\item{[HM2]} Howard, B.J. Millson, J. Snowden, A. Vakil, R.:
{\it The moduli space of n points on the line is
cut out by simple quadrics when\/ $n$ is not six},
ArXiv Mathematics
e-prints math.AG/0607372  (2007)
\medskip
\item{[Ig1]} Igusa, J.:
{\it On the graded ring of theta-constants \/}
Am. J. Math. {\bf 86}, 219-246 (1964).
\medskip
\item{[Ig2]} Igusa, J.:
{\it Modular forms and projective invariants \/}
Am. J. Math. {\bf 89}, 817-855 (1967).
\medskip
\item{[Ka]} Kapranov, M.M.: {\it Chow quotients of Grassmannians I\/},
Adv. Sov. Math. {\bf 16} (2), 29-110 (1993)
\medskip
\item{[Ko1]} Kondo, S.: {\it The moduli space of
Enriques surfaces and Borcherds products,\/}
J. Algebraic Geometry  {\bf 11}, 601-627(2002)
\medskip
\item{[Ko2]} Kondo, S.: {\it  The moduli
space of 8 points on $P^{1}(\cz)$ and automorphic forms\/} in  Algebraic Geometry Contemporary Mathematics {\bf 422}, 89-106
(2007), Amer. Math.Soc.
\medskip
\item{[Koi]} Koike, K.: {\it  The projective
embedding of the configuration space $X(2,8)$\/},
Technical Reports of Mathematical Sciences, Chiba University, {\bf16} (2000)
\medskip
\item{[MY]} Matsumoto, K. Yoshida, M.: {\it Configuration space of 8
points on the projective line and a 5-dimensional Picard modular group,\/}
Compositio Math.~{\bf 86}, 265--280 (1993)
\medskip
\item{[Mu]} {\it Tata Lectures on Theta II\/},
Modern Birkhäuser Classics,
Reprint of the 1984 ed., XIV  (2007)
\medskip
\item{[Ru1]} Runge, B.:
{\it On Siegel modular form, part I\/}, J. Reine angew. Math. {\bf 436}, 57-85 (1993)
\medskip
\item{[Ru2]} Runge, B.:
{\it On Siegel modular forms, part II\/}, Nagoya Math. J. {\bf138}, 179-197 (1995)
\item{[Ts]} Tsuyumune, S.:
{\it Thetanullwerte on a moduli space of curves and hyperelliptic loci\/},   Math. Zeit. {\bf 207}, 539-568 (1991)
\bye